\documentclass[a4paper,12pt]{article}
\pdfoutput=1
\usepackage[titletoc,toc,title]{appendix}
\usepackage[utf8]{inputenc}
\usepackage{amsfonts}
\usepackage{graphicx,color,epstopdf}
\usepackage{amsmath,amsthm, morefloats, bm}
\usepackage[inline]{showlabels}
\allowdisplaybreaks
\newtheorem{prop}{Proposition}[section]

\newtheorem{mytheorem}{Theorem}[section]

\newtheorem{myremark}{Remark}[section]

\def\XXint#1#2#3{{\setbox0=\hbox{$#1{#2#3}{\int}$}
\vcenter{\hbox{$#2#3$}}\kern-.5\wd0}}

\begin{document}
\title{A numerical mode matching method for wave scattering in a layered medium
  with a stratified inhomogeneity}

\author{Wangtao Lu 
	\thanks{School of Mathematical Sciences, Zhejiang University, Hangzhou 310027,
    China.  Email: wangtaolu@zju.edu.cn (corresponding author).} \and
    Ya Yan Lu\thanks{Department of Mathematics, City University of Hong Kong, Kowloon,
      Hong Kong, China. Email: mayylu@cityu.edu.hk.} \and
{Dawei Song}
    \thanks{Department of Mathematics, Nanjing University of Aeronautics and
      Astronautics, Nanjing, Jiangsu, China. Email: dwsmath@nuaa.edu.cn.}
      }
\maketitle
\begin{abstract}
Numerical mode matching (NMM) methods are widely used for analyzing
wave propagation and scattering in structures that are piecewise
uniform along one spatial direction. For open structures that are
unbounded in transverse directions (perpendicular to the uniform
direction), the NMM methods use the perfectly matched layer (PML)
technique to truncate the transverse variables.
When incident waves are specified in homogeneous media 
surrounding the main structure, the total field is not 
always outgoing, and the NMM methods rely on reference
solutions for each uniform segment. Existing NMM methods have difficulty
handing gracing incident waves and special incident waves related to
the onset of total internal reflection, and are not very efficient at
computing reference solutions for non-plane incident waves. In this paper, a new NMM method is developed to overcome these limitations.
A Robin-type boundary condition is proposed to
ensure that non-propagating and non-decaying wave field components are
not reflected by truncated PMLs. Exponential convergence of the PML
solutions based on the hybrid Dirichlet-Robin boundary condition is
established theoretically. A fast method is developed for computing
reference solutions for cylindrical incident waves. The new NMM is implemented for  
two-dimensional structures and polarized electromagnetic waves.  
Numerical experiments are carried out to validate the new NMM method
and to demonstrate  its performance. 
\end{abstract}


\section{Introduction}


Wave scattering problems in a layered medium with a penetrable or
impenetrable inhomogeneity appear in numerous scientific and
engineering applications \cite{chew95}. Classical numerical methods
such as the finite difference method, the finite element method (FEM)
\cite{mon03}, and the spectral method are very versatile, but are not
always the most efficient, since they need to discretize the whole
computational domain.  For piecewise homogeneous structures, the
boundary integral equation (BIE) methods \cite{cai02, brulyoperaratur16,
  laigreone16, luluqia18} are highly competitive since
they discretize only the interfaces and the boundary of the
inhomogeneity.  If the structure can be divided into a number of
segments or regions where the governing equation becomes separable,
the mode matching method, a.k.a mode expansion method or modal
method \cite{botcramcp81, li93a, shestesan82}, and its many numerical variants
\cite{chiyehshi09, gra99, gragui96, kno78, lalmor96, li96, lushilu14, mor95,
  sonyualu11, derdezoly98, biebae01, biederbaeolydez01} may be used.  Typically, these
methods are applicable if the structure is piecewise uniform along one
spatial direction. In each uniform segment, the wave field is expanded
in eigenmodes of a related transverse differential operator, and the
expansion coefficients are solved from a linear system obtained by
matching the wave field at the interfaces between neighboring
segments.  The classical mode matching method solves the eigenmodes
analytically. The numerical mode matching (NMM) methods solve the
eigenmodes by numerical methods, and they are easier to implement and
applicable to more general structures.  The mode matching method and
its variants have the advantage of avoiding discretizing one spatial
variable. They are widely used in engineering applications, since many
designed structures are indeed piecewise uniform. 

For numerical simulations of waves, the perfectly
matched layer (PML) \cite{ber94} is an important technique for
truncating unbounded domains. It is widely used with standard numerical
methods, such as FEM, that discretize the whole computational 
domain.  The BIE methods  usually automatically take care of the radiation
conditions at infinity, but for scattering problems in layerd media,
PML can also be used to efficiently truncate interfaces that extend to
infinity \cite{luluqia18}. For NMM methods, PML was first applied to
study piecewise uniform waveguides \cite{derdezoly98, biebae01, biederbaeolydez01}. An optical waveguide is an open
structure, i.e., the transverse domain perpendicular to the
waveguide axis is unbounded. Analytic mode matching method is
difficult to use, since the transverse operator has a continuous
spectrum and field expansions contain integrals related 
to the radiation modes. When a PML is used to truncate the transverse
domain, typically with a zero Dirichlet boundary condition at the
external boundary of the PML, the continuous spectrum is discretized,
and the field  
expansions are approximated by sums of discrete eigenmodes. 


For many applications, an incident wave is specified in the
homogeneous media surrounding the scatterer, then the total wave field
in each uniform segment does not satisfy outgoing radiation conditions
in the transverse directions, and is incompatible with the eigenmodes
computed using a PML.  To overcome this difficulty, we can find a
reference solution for the given incident wave in each uniform
segment, and then expand the difference between the total field and
the reference solution in the PML-based eigenmodes \cite{lushilu14}.
Typically, the field difference in each segment is indeed outgoing in
the transverse directions, and a NMM method based on this approach
works reasonably well. But unfortunately, the method breaks down in
special circumstances where the field difference in a segment contains
a component that is exactly or nearly invariant in the transverse
direction, i.e., a component with a zero transverse wavenumber.  This
happens if the incident wave has the critical incident angle for the
onset of total internal reflection in the exterior segments. In
that case, the field difference in any interior segment contains a
non-propagating and non-decaying component with a zero or near zero
transverse wavenumber.  This difficulty also arises when the incident
wave is nearly parallel to the uniform direction, i.e., a gracing
incidence.  In that case, the field difference in an interior segment
also contains a plane wave component with a near zero transverse
wavenumber. 

In this paper, we develop a new NMM method to overcome the above
difficulty. Our approach is to use a Robin boundary condition for the
PML in the interior segments. The boundary condition is designed
to ensure that the field component with a zero or near zero transverse
wavenumber is not reflected by the PML. A similar Robin-type condition
for PML was previously used by one  of the authors to preserve a
weakly confined guided mode propagating in optical waveguides
\cite{holu03}. For the exterior segments, we keep the simple zero
Dirichlet boundary condition. To give the method a theoretical 
foundation, we analyze the effectiveness of the PML using hybrid
Dirichlet-Robin boundary conditions. 
It is shown that the error induced by the PML decays  exponentially 
with the thickness or the absorbing coefficient of the PML. 
For scattering problems with incident waves from a point or line
source, the NMM method faces an additional difficulty, namely, the
computation of the reference solutions, especially for the segment
involving the inhomogeneity. The traditional approach that turns a
point or line source to plane waves by Fourier transform is not very
efficient. We develop an efficient method for
computing the reference solutions based the PML technique and the
method of separation of variables. 

The rest of this paper is organized as follows. In Section 2, we formulate the
scattering problem, review the PML theory. In Section 3, we describe an NMM
method. In Section 4, we derive the new Robin-type boundary condition and show
that the solution based on a PML and a hybrid Dirichlet-Robin condition
converges to the true scattering solution exponentially. In Section 5, we
develop an efficient method for computing reference solutions when the incident
wave is a line source. In Section 6, we present a few numerical examples to
validate our NMM method and to illustrate its performance. The paper is
concluded by some remarks and discussions in Section 7.


\section{Problem formulation}
To simplify the presentation, we begin with a scattering problem in a
two-layer medium. The physical structure is characterized by a $z$-invariant
dielectric function
\begin{align}
  \label{eq:eps:func}
  \varepsilon(x,y) = \left\{
  \begin{array}{lc}
    \varepsilon_{+}=n_+^2, & (x,y)\in \mathbb{R}_+^2 \backslash \bar{D},\\
    \tilde{\varepsilon}(y), & (x,y)\in D,\\
    \varepsilon_{-}=n_-^2, & (x,y)\in \mathbb{R}_-^2\backslash \bar{D},
  \end{array}
  \right.
\end{align}
where $\mathbb{R}^2_{\pm}=\{(x,y)\in\mathbb{R}^2:\pm y>0\}$, 
$\tilde{ \varepsilon }(y)\geq 1$ is piecewise smooth on 
$(y_0,y_1)$,  $D$ is a rectangle $(-x_0,x_0) \times (y_0,y_1)$ with $x_0>0$,
$y_1\geq 0$ and $y_0\leq 0$, and it corresponds to a stratified
inhomogeneity.  
In $\mathbb{R}_+^2/\bar{D}$, we specify a plane incident wave $u^{\rm
  inc}=e^{i(\alpha x -\beta_+ y)}$, where $\alpha=k_0n_+\sin\theta$, 
$\beta_+=k_0n_+\cos\theta$, 
and $\theta\in(-\pi/2,\pi/2)$ is the incident angle. The
total wave field $u^{\rm tot}$ satisfies the Helmholtz equation
\begin{align}
  \label{eq:gov:problem}
  \Delta u^{\rm tot} + k_0^2 \varepsilon(x,y) u^{\rm tot} = 0,
\end{align}
where $\Delta=\partial_x^2 + \partial_y^2$, and $k_0$ is the free-space
wavenumber.  Across an interface or discontinuity, we have the
following transmission condition
\begin{align}
  \label{eq:trans:cond}
  [u^{\rm tot}] = 0, \quad \left[\frac{\partial u^{\rm tot}}{\partial{\bm \nu}} \right]= 0,
\end{align}
where ${\bm \nu}$ is the unit normal vector on the interface pointing towards
$\mathbb{R}_-^2$, and $[\cdot]$ denotes the jump of the quantity across the
interface. For electromagnetic waves in the $E$ polarization,
  $u^{\rm 
  tot}$ is the $z$ component (the only nonzero component) of the 
electric field.

At infinity, the scattered wave field  defined as
\[
  u^s=\left\{ \begin{array}{lc}
                u^{\rm tot} - u^{\rm tot}_0, & {\rm in}\ \mathbb{R}^2/\bar{D},\\
                u^{\rm tot}, &{\rm in}\  D,
                \end{array}
  \right.
\]
satisfies the half-plane Sommerfeld radiation condition in both $\mathbb{R}^2_+$
and $\mathbb{R}^2_-$, i.e.,
\begin{align}
  \label{eq:half:Som:cond}
\lim_{r\rightarrow \infty}\sqrt{r}(\partial_r u^s - ik_0n_{\pm}u^s) = 0,\quad r=\sqrt{x^2+y^2},\quad (x,y)\in\mathbb{R}_{\pm}^2/\bar{D}.
\end{align}
Here, $u^{\rm tot}_0$ is the solution for the same incident wave in
the background two-layer medium without the inhomogeneity. More
precisely, we have
\begin{align}
  \label{eq:sol:cond:c}
  u^{\rm tot}_0 &= \left\{\begin{array}{lc}
                           e^{i(\alpha x - \beta_+ y)} + Re^{i(\alpha x + \beta_+ y)} &(x,y)\in\mathbb{R}_+^2,\\
                           (R+1)e^{i(\alpha x -  \beta_-y)} &(x,y)\in\mathbb{R}_-^2,
                         \end{array}\right.\\
  \label{eq:sol:para}
  \beta_- &= \sqrt{k_0^2\varepsilon_--\alpha^2},\quad R =
        \frac{\beta_+-\beta_-}{\beta_++\beta_-}.
\end{align}
According to \cite{roazha92,chezhe10,baohuyin18}, we have the following existence and uniqueness results:
\begin{mytheorem}
  For any incident plane wave with $k_0>0$,
  the scattering  problem (\ref{eq:gov:problem}), (\ref{eq:trans:cond}),
  (\ref{eq:half:Som:cond}) has a unique solution $u^{\rm tot}$ 
in $H_{\rm loc}^{1}(\mathbb{R}^2)$.
\end{mytheorem}

Since $u^s$ is outgoing, the PML technique \cite{ber94} can be used to
truncate $\mathbb{R}^2$. Let us define 
the following complex coordinate stretching functions
\begin{align}
  \label{eq:pml:xy}
  \tilde{x}(x) = x + i\int_{0}^{x}\sigma_1 (t) dt,\quad
  \tilde{y}(y) = y + i\int_{0}^{y}\sigma_2 (t) dt,
\end{align}
where
$\sigma_l(t) = \sigma_l(-t)$ for all $t$, $\sigma_l(t)=0$ for $|t|\leq
L_l/2$, and $\sigma_l(t)>0$ for $|t|>L_l/2$, 
and $L_l>0$ for $l=1$, $2$. Notice that 
the rectangle $B_1= (-L_1/2,L_1/2) \times (-L_2/2,L_2/2)$ encloses the
inhomogeneity $D$, and the rectangle $B_2=(-L_1/2-d_1,L_1/2+d_1)
\times (-L_2/2-d_2,L_2/2+d_2)$ is used to truncate
$\mathbb{R}^2$.  Based on Green's representation formula, the
extension of $u^s$ in $B_2$ can be defined, then 
\[
  \tilde{u}^s(x,y):=u^s(\tilde{x}(x),\tilde{y}(y))
\]
satisfies the following PML-Helmholtz equation
\begin{align}
  \label{eq:upml:1}
  &\nabla\cdot(A\nabla\tilde{u}^s) + \alpha_1(x)\alpha_2(y)k_0^2\varepsilon(x,y)\tilde{u}^s = 0,\\
  \label{eq:upml:2}
  &\tilde{u}^{s}(x,0+) = \tilde{u}^{s}(x,0-), \quad \partial_{\tilde{y}} \tilde{u}^s(x,0+) = \partial_{ \tilde{y} } \tilde{u}^s(x,0-),\ {\rm in}\ \mathbb{R}/[-x_0,x_0].
\end{align}
where $A={\rm diag}(\alpha_2(y)/\alpha_1(x),\alpha_1(x)/\alpha_2(y))$, and
$\alpha_l=1+i\sigma_l$. Typically, a zero Dirichlet boundary condition is
enforced on $\Gamma_2=\partial B_2$, i.e.,
\begin{align}
  \label{eq:upml:3}
  \tilde{u}^s(x,y) = 0,\quad{\rm on}\quad \Gamma_2.
\end{align}
The following theorem characterizes the exponential convergence of the PML
solution.
\begin{mytheorem}
  Let $\sigma_1=\sigma_2\equiv\sigma$ and $d_1=d_2 =d$ in the PML for some
  positive constants $\sigma$ and $d$, and let $\bar{\sigma}=\sigma d$ such that
  $\gamma_0\bar{\sigma}\geq \max(k_{\rm min}^{-1}, d)$, where $\gamma_0 = d /
  \sqrt{(L_1+d)^2+(L_2+d)^2}$. Then for sufficiently large $\bar{\sigma}$, the PML
  problem (\ref{eq:upml:1}-\ref{eq:upml:3}) has a unique solution $\tilde{u}^s$
  in $H^1(B_2)$. Moreover, there exists a constant $C$, which depends only on
  $\gamma_0$, $k_{\rm max}/k_{\rm min}$, and $L_2/L_1$, but independent of $n_-$, $n_+$, $L_1$,
  $L_2$, and $d$, such that
  \begin{align}
    \label{eq:est:ut:u}
    ||u^s-\tilde{u}^s||_{H^1(B_1)}\leq C(1+\hat{C}^{-1})\gamma_1(1+k_{\rm min}L_1)^3\alpha_m^3\nonumber\\
    (1+\bar{\sigma}/L_1)^2e^{-k_{\rm min}\gamma_0\bar{ \sigma }}||\tilde{u}^s||_{H^{1/2}(\Gamma_1)},
  \end{align}
  where $k_{\rm min}=k_0\min(n_-,n_+)$, $k_{\rm max}=k_0\max(n_-,n_+)$,
  $\gamma_1=e^{L_2\sqrt{|k_{\rm max}^2-k_{\rm min}^2|}/2}$,
  $\alpha_m=\sqrt{1+\sigma^2}$, and $\hat{C} =
  \frac{\min(1,\sigma^3)}{2(1+\sigma^2)^2\max(1,k_{\rm max}^2d^2)}$.
  \begin{proof} The PML problem can be considered in $B_2/\bar{D}$ by
    regarding $D$ as an obstacle and enforcing the Dirichlet
    boundary condition $\tilde{u}^s|_{\partial D}=u^s|_{\partial
      D}$. Evidently, this theorem follows directly from Theorem 7.2 in
    \cite{chezhe10}. 
  \end{proof}
\end{mytheorem}
 
Thanks to Theorem 2.2, $\tilde{u}^s$ converges to $u^s$ exponentially in
$B_1$. Therefore, we only need to deal with $\tilde{u}^s$ in the
bounded domain $B_2$ instead of $u^s$ in $\mathbb{R}^2$.

\section{Numerical mode matching method}

For the scattering problem formulated above, the NMM methods are
applicable, since the structure is uniform in $x$ in three different
segments corresponding to $x<-x_0$, $-x_0 < x < x_0$ and $x> x_0$,
respectively. Since a PML is used in the NMM method, we 
define the three segments by $S_1=\{(x,y)|-d_1-L_1/2<x< -x_0\}\cap 
B_2$, $S_2=\{(x,y)|-x_0< x< x_0\}\cap B_2$ and 
$S_3=\{(x,y)|x_0<x<d_1+L_1/2\}\cap B_2$. It is clear that
$\varepsilon(x,y)=\varepsilon_i(y)$ in $S_i$ is independent of $x$,
for $i=1$, 2, 3. 
Accordingly, $\Gamma_2$ (the boundary of $B_2$) is decomposed into 
three parts $\Gamma_{2}^1$, $\Gamma_{2}^2$ and $\Gamma_{2}^3$. In 
particular, 
\[
  \Gamma_2^2=\{(x,d_2+L/2)||x|<x_0\}\cup\{(x,-d_2-L/2)||x|<x_0\}. 
\]

In the last several decades, many different NMM methods have been
developed. These methods use different numerical methods to solve the
eigenmodes in the uniform segments, and also use different techniques
to impose the continuity conditions at the interfaces between the
neighboring segments. Our NMM method is similar to the one presented in
\cite{lushilu14}, and its basic steps are summarized below.

We consider segments $S_1$ and $S_3$ first. According to
Eqs.~(\ref{eq:upml:1}-\ref{eq:upml:3}), $\tilde{u}^s$ in $S_i$ ($i=1$,
3) solves 
\begin{align}
  \label{eq:us:i1}
  &\nabla\cdot(A\nabla\tilde{u}^s) + k_0^2\alpha_1(x)\alpha_2(y)\varepsilon_i(y) \tilde{u}^s = 0,\\
  \label{eq:us:i2}
  &\tilde{u}^{s}(x,0+) = \tilde{u}^{s}(x,0-), \quad \partial_{\tilde{y}} \tilde{u}^s(x,0+) = \partial_{ \tilde{y} } \tilde{u}^s(x,0-),\\
  \label{eq:us:i3}
  &\tilde{u}^{s}(x, d_2+L_2/2) = \tilde{u}^{s}(x,-d_2-L_2/2)=0.
\end{align}
By the method of separation of variables, inserting 
$\tilde{u}^s(x,y)=\phi(y)\psi(x)$ into
(\ref{eq:us:i1}-\ref{eq:us:i3}), we obtain the following eigenvalue
problem for $\phi(y)$ 
\begin{align}
  \label{eq:phi:te}
  &\frac{1}{\alpha_2}\frac{d}{dy}\left( \frac{1}{\alpha_2}\frac{d\phi}{dy} \right) + k_0^2\varepsilon_i(y)\phi(y)=\delta\phi,\\
  \label{eq:phi:te2}
  &\phi(0+) = \phi(0-), \phi'(0+) = \phi'(0-),\\
  \label{eq:phi:te3}
  &\phi(d_2+L_2/2) = \phi(-d_2-L_2/2) = 0,
\end{align}
and the associated equation for $\psi(x)$
\begin{align}
  \label{eq:psi:te}
  \frac{1}{\alpha_1}\frac{d}{dx}\left( \frac{1}{\alpha_1}\frac{d\psi}{dx} \right) + \delta\psi = 0.
\end{align}

The above Sturm-Liouville eigenvalue problem (\ref{eq:phi:te}-\ref{eq:phi:te3}) for $\phi$ is not
self-adjoint, thus $\delta$ is in general complex. Nevertheless,
$\delta$ can be forced to the upper half-plane based on the following
proposition. 
\begin{prop}
  Under the same assumptions as Theorem 2.2, we have that for sufficiently large
  $\bar{\sigma}$, ${\rm Im}(\delta)\geq 0$.
  \begin{proof}
    See Proposition A.1 in Appendix A.
  \end{proof}
\end{prop}
As in \cite{sonyualu11}, we employ a pseudospectral method 
 to find the numerical eigenmodes. 
Assuming $N$ eigenpairs $\{\delta_j, \phi_j(y)\}$ for $j=1$, ..., $N$,
are obtained based on the $N$ collocation points
$\{y^j\}_{j=1}^{N}\subset [-d_2-L_2/2, d_2+L_2/2]$, we
approximate 
$\tilde{u}^s$ by
\begin{align}
  \label{eq:region1}
  \tilde{u}^s \approx \sum_{j=1}^{N}\left[ c_j^{(1)}
  e^{-i\sqrt{\delta_j}(\tilde{x}(x)-\tilde{x}(-x_0))} + d_j^{(1)}
  e^{i\sqrt{\delta_j}(\tilde{x}(x)-\tilde{x}(-d_1-L_1/2))} \right] \phi_j(y)
\end{align}
in $S_1$, and by
\begin{align}
  \label{eq:region3}
  \tilde{u}^s \approx \sum_{j=1}^{N} \left[ c_j^{(3)}
  e^{-i\sqrt{\delta_j}(\tilde{x}(x)-\tilde{x}(d_1+L_1/2))} + d_j^{(3)}
  e^{i\sqrt{\delta_j}(\tilde{x}(x)-\tilde{x}(x_0))} \right] \phi_j(y)
\end{align}
in $S_3$, where $\sqrt{\delta_j}$ is defined to be in the branch with ${\rm
  Im}(\sqrt{\delta_j})\geq 0$ and hence with ${\rm Re}(\sqrt{\delta_j})\geq 0$
according to Proposition 1. Based on the zero Dirichlet boundary
condition at 
$x=\pm(d_1+L_1/2)$, we get
\begin{align}
  \label{eq:dj1}
 d_j^{(1)} &= -c_j^{(1)} e^{-i\sqrt{\delta_j}(\tilde{x}(-d_1-L_1/2)-\tilde{x}(-x_0))},\\
  \label{eq:cj3}
 c_j^{(3)} &= -d_j^{(3)}e^{i\sqrt{\delta_j}(\tilde{x}(d_1+L_1/2)-\tilde{x}(-x_0))}.
\end{align}
Therefore, 
\begin{align*}
  |d_j^{(1)}|=|c_j^{(1)}| e^{-{\rm Im}(\sqrt{ \delta_j })(L_1/2-x_0+d_1) - {\rm Re}(\sqrt{ \delta_j })\int_{L_1/2}^{L_1/2+d_1}\sigma(t)dt}\approx 0,\\
  |c_j^{(3)}|=|d_j^{(3)}| e^{-{\rm Im}(\sqrt{ \delta_j })(L_1/2-x_0+d_1) - {\rm Re}(\sqrt{ \delta_j })\int_{L_1/2}^{L_1/2+d_1}\sigma(t)dt}\approx 0,
\end{align*}
for sufficiently large $\sigma$ and $d_1$. 
Consequently, we can assume that there are no terms with coefficients
$d_j^{(1)}$ and  $c_j^{(3)}$ in Eqs.~(\ref{eq:region1}) and
(\ref{eq:region3}), respectively. Physically, this corresponds to the
fact that $u^s$ should not contain incoming waves in the two
exterior segments. 

In segment $S_2$, we have
\begin{align}
  \label{eq:eps2}
  \varepsilon_2(y) = \left\{
    \begin{array}{lc}
    n_+^2, & y>y_1,\\
    \tilde{\varepsilon}(y), & y_0<y<y_1,\\
    n_-^2, & y<y_0.
    \end{array}
  \right.
\end{align}
The method of separation of variables is not applicable to
$\tilde{u}^s$, since it does not satisfy the homogeneous transmission
conditions (\ref{eq:trans:cond}) at $y=y_0$ and $y=y_1$. Instead, we 
need to subtract from $u^{\rm tot}$ a wave field $u^{\rm tot}_2$
which solves the scattering problem for the same incident wave and a
layered profile $\varepsilon(x,y)=\varepsilon_2(y)$ in
$\mathbb{R}^2$. 
We let $u^{\rm tot}_2$ be the solution with the same $x$-dependence as
the incident wave. More details are given Proposition A.2.

For $u_2^s = u^{\rm tot} - u_2^{\rm tot}$, we enforce the same zero
Dirichlet boundary condition 
\begin{equation}
  \label{eq:dir:cond}
  \tilde{u}_2^s=0,\quad{\rm on}\ \Gamma_2^2,
\end{equation}
where $\tilde{u}_2^s(x,y)=u_2^s(\tilde{x}(x),\tilde{y}(y))$. The method of
separation of variables can be applied to $\tilde{u}_2^s$. 
Based on the same discretization points $\{y^j\}_{j=1}^N$, 
we obtain $N$ eigenpairs
$\{\delta_j^{(2)},\phi_j^{(2)}(y)\}_{j=1}^{N}$ in $S_2$, then 
\begin{align}
  \label{eq:region2}
  \tilde{u}_2^s \approx \sum_{j=1}^{N}\left[ c_j^{(2)} e^{-i\sqrt{\delta_j^{(2)}}(x-x_0)} + d_j^{(2)}
  e^{i\sqrt{\delta_j^{(2)}}(x+x_0)}  \right] \phi_j^{(2)}(y).
\end{align}

On the two interfaces between $S_2$ and the other two segments $S_1$
and $S_3$, i.e.  at $x=\pm x_0$, we have the transmission conditions
\begin{align}
  \label{eq:gov:tran1}
  [\tilde{u}_s(\pm x_0,y) - \tilde{u}_s^2(\pm x_0,y)] &= [\tilde{f}(\pm x_0,y)],\\
  \label{eq:gov:tran2}
  [\partial_x \tilde{u}_s(\pm x_0,y) - \partial_x \tilde{u}_s^2(\pm
  x_0,y)] &= [\tilde{g}(\pm x_0,y)], 
\end{align}
where $\tilde{f}(x,y)=f(x,\tilde{y}(y))$, $\tilde{g}(x,y)=g(x,\tilde{y}(y))$, and
\[
  f(x,y)=u^{\rm tot}_2(x,y) - u^{\rm tot}_0(x,y),\quad  g(x, y)=\partial_x u^{\rm tot}_2(x,y) - \partial_x u^{\rm tot}_0(x,y).
\]
Collocating (\ref{eq:gov:tran1}) and (\ref{eq:gov:tran2}) at $y=y^j$ for
$j=1,\ldots, N$, and using Eqs.~(\ref{eq:region1}), (\ref{eq:region3}) and
(\ref{eq:region2}), we obtain a linear system
\begin{align}
  {\bm A}\left[\begin{array}{c}
                 {\bm c}^{(1)}\\
                 {\bm c}^{(2)}\\
                 {\bm d}^{(2)}\\
                 {\bm d}^{(3)}\\
                 \end{array}\right] = {\bm b},
\end{align}
where ${\bm A}$ is a $4N\times 4N$ matrix, ${\bm b}$ is a
$4N\times 1$ matrix, ${\bm c}^{(i)}=[c_1^{(i)}, \ldots,
c_N^{(i)}]^{T}$, etc. 
Solving the above system, we get $\tilde{u}^s$ in $S_1$ and $S_3$, and
$\tilde{u}^s_2$ in $S_2$, thus $u^{\rm tot}$ can be found in
the physical domian $B_1$. 

In the above, the NMM method is only presented for the case of a single
inhomogeneous segment in a two-media layered background. It is straightforward
to extend the NMM method to structures with multiple inhomogeneous segments that
are uniform along the same direction. The NMM method can also be used to study
scattering problems in the $H$ polarization (the only nonzero component of the
magnetic field is its $z$ component) and problems involving perfect electrical
conductor (PEC) or perfect magnetic conductor (PMC) scatterers.

\section{A Robin-type boundary condition}
As we mentioned in the introduction, the NMM method based on the zero
Dirichlet condition (\ref{eq:dir:cond}) usualy works,  but in some
special circumstances, it  exhibits a slow convergence and even a
divergence, since $u_2^s$ may not be strictly outgoing. It should be
pointed out that there is no contradiction with Theorem 2.2, since
that theorem is about applying the PML to $u^s$, but the NMM method
applies the PML to $u_2^s$ for the interior segment $S_2$. 

In fact, it is easy to deduce that
\begin{align}
  \label{eq:us2:asym:p}
  u_2^s&= (R-R_2e^{-2i\beta_+ y_1})e^{i(\alpha x+\beta_+ y)} + u^s,\quad{\rm for}\
  y>y_1,\\
  u_2^s&=(T - T_2e^{i\beta_-y_0-i\beta_+ y_1})e^{i(\alpha x - \beta_-
         y)}+u^s,\quad{\rm for}\ y<y_0, 
\end{align}
where $R$ and $R_2$ ($T$ and $T_2$) are the reflection (transmission)
coefficients in the reference solutions $u_0^{\rm tot}$ and $u_2^{\rm
  tot}$, respectively. Therefore, $u_2^s$ can be decomposed as a
scattered wave field and an up-going plane wave with $y$-wavenumber
$\beta_+$ for $y>y_1$ or a down-going plane wave with $y$-wavenumber
$\beta_-$ for $y<y_0$. Consequently, only when $\beta_+$ and 
$\beta_-$ are sufficiently far away from zero, does $u_2^s$ attenuate in the
PML. However, this is not ture for the following two cases:
\begin{itemize}
\item[a.] For gracing incidences with $\theta$ close to
  $\pm\pi/2$,   $\beta_+$ is close to $0$;
\item[b.] For $n_-< n_+$ and at the critical angles $\theta = \pm
  \arcsin (n_-/n_+)$ for the onset of total internal reflection,  $\beta_-=0$.
\end{itemize}
Notice that $\tilde{u}^s\approx 0$ at the exterior boundary of the
PML, therefore
\begin{align*}
  \tilde{u}_2^s&\approx (R-R_2e^{-2i\beta_+ y_1})e^{i(\alpha x+\beta_+
    \tilde{y})},\quad{\rm on}\ y=d_2+L_2/2,\\
  \tilde{u}_2^s&\approx(T - T_2e^{i\beta_-y_0-i\beta_+ y_1})e^{i(\alpha x -
    \beta_- \tilde{y})},\quad{\rm on}\ y=-d_2-L_2/2.
\end{align*}
However, the NMM method is not compatible with the above
inhomogeneous boundary conditions. To overcome this difficulty, the
following result is needed. 
\begin{prop}
The scattered fields $u^s$ and $u_2^s$ satisfy 
 \begin{align}
   \label{eq:us:u}
   \partial_{y} u_2^s - i\beta_+ u_2^s = \partial_{y} u^s - i\beta_+ u^s,\quad{\rm on}\ y=d_2+L_2/2,\\
  \label{eq:us:b}
   \partial_{y} u_2^s + i\beta_- u_2^s = \partial_{y} u^s + i\beta_- u^s,\quad{\rm on}\ y=-d_2-L_2/2.
\end{align}
Correspondingly,
\begin{align}
  \label{eq:tus:u}
   (A\nabla \tilde{u}_2^s)\cdot{\bf \nu} - i\beta_+ \tilde{u}_2^s =
(A\nabla\tilde{u}^s)\cdot{\bf \nu} - i\beta_+ \tilde{u}^s,\quad{\rm on}\
y=d_2+L_2/2,\\
  \label{eq:tus:b}
   (A\nabla \tilde{u}_2^s)\cdot{\bf \nu} - i\beta_- \tilde{u}_2^s =
(A\nabla\tilde{u}^s)\cdot{\bf \nu} - i\beta_- \tilde{u}^s,\quad{\rm on}\
y=-d_2-L_2/2.
 \end{align}
 \begin{proof}
   The proof is straightforward.
 \end{proof}
\end{prop}
The above proposition suggests the following homogeneous Robin boundary
conditions,
\begin{align}
  \label{eq:cond:tus:1}
  \frac{1}{\alpha_2}\partial_y\tilde{u}_2^s-i\beta_+\tilde{u}_2^s&=(A\nabla \tilde{u}_2^s)\cdot{\bf \nu} - i\beta_+ \tilde{u}_2^s = 0,\quad{\rm on}\quad y=d_2+L_2/2,\\
  \label{eq:cond:tus:2}
   \frac{1}{\alpha_2}\partial_y\tilde{u}_2^s+i\beta_+\tilde{u}_2^s&=(A\nabla \tilde{u}_2^s)\cdot{\bf \nu} - i\beta_- \tilde{u}_2^s = 0,\quad{\rm on}\quad y=-d_2-L_2/2.
\end{align}
Based on the pseudospectral method \cite{sonyualu11}
and the above boundary conditions, we can find the eigenmodes
$\phi_j^{(2)}$ ($1 \le j \le N$), and expand  $\tilde{u}_2^s$  in $S_2$
in these eigenmodes.

Although different boundary conditions are used on $\Gamma_2$, 
the following theorem ensures that  $\tilde{u}^s$ still converges 
to $u^s$  exponentially. 
\begin{mytheorem}
  Under the same assumptions as Theorem 2.2, we have that for sufficiently large
  $\bar{\sigma}$, the PML problem (\ref{eq:upml:1}), (\ref{eq:upml:2}) equipped
  with the following hybrid Dirichlet-Robin boundary condition
  \begin{align}
    \label{eq:general:cond}
    \left\{\begin{array}{lc}
             \tilde{u}^s = 0,&{\rm on}\quad \Gamma_2/\overline{{\Gamma}}\\
             ( A\nabla\tilde{u}^s )\cdot {\bm \nu} = iW^2 \tilde{u}^s,&{\rm on}\quad \Gamma,
           \end{array}
                                                                            \right.
  \end{align}
  where $W\in L^{\infty}(\Gamma)$ is real-valued and
  $\Gamma\subset\Gamma_2$ is an open bounded set, has a unique
  solution $\tilde{u}^s$ in $H^1(B_2)$. Moreover, there exists a constant $C$,
  which depends only on $||W^2||_{L^\infty(\Gamma)}$, $\gamma_0$, $k_{\rm
    max}/k_{\rm min}$, and $L_2/L_1$, but independent of $n_-$, $n_+$, $L_1$,
  $L_2$, and $d$, such that
  \begin{align}
    ||u^s-\tilde{u}^s||_{H^1(B_1)}\leq C(1+\hat{C}^{-1})\gamma_1(1+k_{\rm min}L_1)^3\alpha_m^3\nonumber\\
    (1+\bar{\sigma}/L_1)^2e^{-k_{\rm min}\gamma_0\bar{ \sigma }}||\tilde{u}^s||_{H^{1/2}(\Gamma_1)}.
  \end{align}
   \begin{proof}
    This can be proved by the similar argument as the proof of Theorem 7.2 in
    \cite{chezhe10}. We here only mention significant modifications. For consistency
    and simplicity, we will load the whole notations from \cite{chezhe10} and will use
    them only in this proof so that $x=(x_1,x_2)$ now denotes a point but not a
    scalar, etc..

    The PML equation in the PML layer (see Eqs. (5.1-5.3) in \cite{chezhe10}) for the
    generalized boundary condition (\ref{eq:general:cond}) should be revised to
    \begin{align}
      \label{eq:proof:pmleq}
      &\nabla\cdot(A\nabla w) + \alpha_1\alpha_2k^2w = 0,\quad {\rm in}\quad \Omega_{\rm PML} = B_2\backslash\bar{B}_1,\\
      \label{eq:proof:pmleq2}
      &[w]_\Sigma=\left[ \frac{\partial w}{\partial x_2} \right]_\Sigma=0\quad {\rm on}\quad \Sigma\cap\Omega_{\rm PML},\\
      \label{eq:proof:pmleq3}
      &w=0\ {\rm on}\ \Gamma_1,\quad w=q\ {\rm on}\ \Gamma_2\backslash\bar{\Gamma},\quad (A\nabla w)\cdot{\bf \nu}-iW^2 w=\tilde{q},\ {\rm on}\quad\Gamma,
    \end{align}
    where $q\in H^{1/2}(\Gamma_2/\bar{\Gamma})$ and $\tilde{q}\in
    H^{-1/2}(\Gamma)$. Then, the related sesquilinear form $c: H^{1}(\Omega_{\rm
      PML})\times H^{1}(\Omega_{\rm PML})\rightarrow \mathbb{C}$ becomes
    \[
      c(\varphi,\psi) =
      \int_{\Omega_{PML}}(A\nabla\varphi\cdot\nabla\bar{\psi}-\alpha_1\alpha_2k^2\varphi\bar{\psi})dx
      - i\int_{\Gamma}W^2\varphi\bar{\psi}ds,\forall\varphi,\psi\in
      H_{0/\Gamma}^1(\Omega_{\rm PML}),
    \]
    where $H_{0/\Gamma}^1(\Omega_{\rm PML}):=\{v\in H^{1}(\Omega_{\rm PML}):
    v=0\ {\rm on}\ \Gamma_1\cup\Gamma_2/\bar{\Gamma}\}$. The weak formulation of
    (\ref{eq:proof:pmleq}-\ref{eq:proof:pmleq3}) is: Find $w\in
    H^{1}(\Omega_{\rm PML})$ such that Eq.(\ref{eq:proof:pmleq3}) is satisfied
    and that
    \begin{align}
      \label{eq:proof:weakform}
      c(w,v) = \langle \tilde{q}, v\rangle|_{\Gamma},\quad\forall v\in H_{0/\Gamma}^1(\Omega_{\rm PML}).
    \end{align}
    Correspondingly, the weighted $H^1$-norm is revised to
    \[
      |||\varphi|||_{H^1(\Omega)}=\left( ||\nabla\varphi||^2_{L^2(\Omega)} +
        ||k\varphi||^2_{L^2(\Omega)} +
        ||W\varphi||_{\Gamma}^2/(1+\sigma^2)^2 \right)^{1/2},
    \]
    and the equivalent norm on $H^{1}(\Omega_{\rm PML})$ becomes
    \[
      ||\varphi||_{*,\Omega_{\rm PML}} = \left(
        ||A\nabla\varphi||^2_{L^2(\Omega_{\rm PML})} +
        ||k\alpha_1\alpha_2\varphi||^2_{L^2(\Omega_{\rm PML})} +
        ||W\varphi||_{\Gamma}^2 \right)^{1/2}.
    \]

    Next, one sees that Lemma 5.1 in \cite{chezhe10} still holds with the space $H_0^1$
    replaced by $H_{0/\Gamma}^1$. The proof relies on the following estimates
    \[
      ||\varphi||_{L^2(\Omega_1)}^2\leq d_1^2\left| \left|
          \frac{\partial\varphi}{\partial x_1} \right|
      \right|^2_{L^2(\Omega_1)},\quad ||\varphi||_{L^2(\Omega_2)}\leq
      d_2^2\left| \left| \frac{\partial\varphi}{\partial x_2}\right| \right|^2,
    \]
    which were proved by using $\varphi=0$ on $\Gamma_2$ in \cite{chezhe10}. However, we
    remark that these two estimates still hold even when $\varphi\neq0$ on
    $\Gamma_2$ since we always have $\varphi=0$ on $\Gamma_1$ for $\varphi\in
    H_{0/\Gamma}^1(\Omega_{\rm PML})$.

    Thus, it is clear that Lemma 5.2 in \cite{chezhe10} holds with the space of $\zeta$
    replaced by ``for any $\zeta\in H^1(\Omega_{\rm PML})$ such that $\zeta=0$
    on $\Gamma_1$, $\zeta=q$ on $\Gamma_2/\bar{\Gamma}$, and
    $A(\nabla\zeta)\cdot{\bf \nu}-iW^2\zeta=\tilde{q}$ on $\Gamma$.''

    Next, Lemma 7.1 in \cite{chezhe10} holds after one replaces $X(f)$ with
    \begin{align*}
      \tilde{X}(f):=\{&\zeta\in H^1(\Omega_{\rm PML}): \zeta=0\ {\rm on}\
                        \Gamma_1, \zeta
                        = \mathbb{E}(f)\ {\rm on}\ \Gamma_2/\bar{\Gamma}, \\
                      &(A\nabla\zeta)\cdot {\bf \nu} - iW^2\zeta =
                        (A\nabla_x\mathbb{E}(f))\cdot {\bf \nu} -
                        iW^2\mathbb{E}(f)\ {\rm on}\ \Gamma\}.
    \end{align*}
    Here, we will have
    \begin{align*}
      \inf_{\zeta\in \tilde{X}}||\zeta||_{*,\Omega_{\rm PML}}\leq
      C(1+k_1L_1)\alpha_m^2(&  ||\mathbb{E}(f)||_{H^{1/2}(\Gamma_2/\bar{\Gamma})} \nonumber\\
                            &+ ||(A\nabla_x\mathbb{E}(f))\cdot{\bf \nu}-iW^2\mathbb{E}(f)||_{H^{-1/2}}(\Gamma)),
    \end{align*}
    where $C$ now depends on the norm $||W^2||_{L^2(\Gamma_2)}$ considering
    the modified norm $||\cdot||_{*,\Omega_{\rm PML}}$. Since $\mathbb{E}(f)$ is
    smooth on $\Gamma$,
    \[
      (A\nabla_x\mathbb{E}(f))\cdot{\bf \nu}-iW^2\mathbb{E}(f)\in
      L^{2}(\Gamma)\cap L^{\infty}(\Gamma),
    \]
    so that
    \[
      ||(A\nabla\mathbb{E}(f))\cdot{\bf
        \nu}-iW^2\mathbb{E}(f)||_{H^{-1/2}(\Gamma)}\leq
      C||\mathbb{E}(f)||_{W^{1,\infty}(\Gamma)},
    \]
    and hence Lemma 7.1 in \cite{chezhe10} follows which proves the theorem.
  \end{proof}
\end{mytheorem}
If we define $W$ in $\Gamma_2^2$ by 
\begin{align*}
  W(x,y) = \left\{
  \begin{array}{lc}
    \sqrt{\beta_+}, & {\rm on}\ y=d_2 + L_2/2,\\
    \sqrt{\beta_-}, & {\rm on}\ y=-d_2 - L_2/2,
  \end{array}
                      \right.
\end{align*}
then Theorem 3 is applicable to our scattering problem.
Consequently, with the hybrid Dirichlet-Robin
condition (\ref{eq:general:cond}), $\tilde{u}^s$ still exponentially
converges to $u^s$ in the physical domain $B_1$.

Theorems 2.2 and 3.1 are established for PMLs with constant and 
equal $\sigma_1$ and $\sigma_2$. In practice, we may set 
$\sigma_1(x)$ and $\sigma_2(y)$ as continuous functions to increase 
flexibility. For example, we may choose 
\begin{align}
  \label{eq:pml:setup}
  \sigma_l(t) = \sigma\left( \frac{t - L_l/2}{d_l} \right)^m,\quad{\rm in}\ L_l/2<|t|<L_l/2+d_l,
\end{align}
for a positive constant $\sigma$ and an integer $m\geq 0$, where $m=0$
corresponds to the constant case.

\section{Cylindrical incident waves}

The NMM methods are typically implemented for plane incident
waves. For other incident waves, such as 
point or line sources and Gaussian beams, the NMM methods may be used
with a Fourier transform that rewrites the incident wave as a
superposition of plane waves. This approach is not very efficient,
since it is necessary to solve the problem for many different incident
plane waves. Alternatively, we can try to find a reference
solution for the given non-plane incident wave in each uniform
segment. This task is nontrivial for the interior segment
corresponding to the inhomogeneity. In the following, we present an
efficient method for computing the reference solutions when the
incident wave is a cylindrical wave generated by a line source. 

The incident cylindrical wave is $u^{\rm
  inc}=\frac{i}{4}H_0^{(1)}(k_0n_+\rho(x,y))$ corresponding to a line
source at $(x^*,y^*)\in\mathbb{R}_+^2/\bar{D}$, 
where $\rho(x,y)=\sqrt{(x-x^*)^2 +
  (y-y^*)^2}$. The governing Helmholtz equation becomes
\begin{align}
  \label{eq:gov:ps}
  \Delta u^{\rm tot} + k_0^2\varepsilon(x,y)u^{\rm tot} = -\delta(x-x^*)\delta(y-y^*).
\end{align}
Considering  the location of the source, we have the following three cases: 
\begin{itemize}
\item[(a)] If $|x^*|<x_0$, we set $u_0^{\rm tot}\equiv 0$ and find a
  nonzero $u_2^{\rm tot}$; 
\item[(b)]  If $|x^*|>x_0$ and $y^*>y_1$, we set $u_2^{\rm tot}\equiv
  0$ and find a nonzero $u_0^{\rm tot}$; 
\item[(c)] If $|x^*|=x_0$, then we have to find nonzero
  $u_0^{\rm tot}$ and $u_2^{\rm tot}$.
\end{itemize}
We consider the typical case (a), where $u^{\rm tot}_2$ must be
computed in segment $S_2$. The NMM method requires $u_2^{\rm tot}$ and
its $x$-derivative at $x=\pm x_0$ to evaluate $\tilde{f}$ and $\tilde{g}$
in Eqs.~(\ref{eq:gov:tran1}) and (\ref{eq:gov:tran2}).

Following the one-dimensional profile $\varepsilon_2(y)$ given in
(\ref{eq:eps2}), $\mathbb{R}^2$ can be split into 
three layers $y<y_0$, $y_0<y<y_1$, and $y>y_1$. The wave field
\[
  u^{s} = \left\{
    \begin{array}{lc}
      u_2^{\rm tot} - u^{\rm inc},&{\rm in}\ y>y_1,\\
      u_2^{\rm tot},&{\rm otherwise},
    \end{array}
  \right.
\]
is outgoing as $y \to \pm \infty$. Using the same PML
as before and applying the technique of separation of variables to $\tilde{u}^s$, we
obtain the following eigenvalue problem
\begin{align}
  \label{eq:eigprob:x1}
  &\frac{1}{\alpha_1}\frac{d}{dx}\left( \frac{1}{\alpha_1} \frac{d\psi}{dx}\right) = \delta \psi,\\
  \label{eq:eigprob:x2}
  &\psi(-L_1/2-d_1) = \psi(L_1/2+d_1) = 0,
\end{align}
and its associated equation
\begin{align}
  \label{eq:asso:phi}
  \frac{1}{\alpha_2}\frac{d}{dy}\left( \frac{1}{\alpha_2} \frac{d\phi}{dy}\right) + (k_0^2\varepsilon_2(y)+\delta) \phi = 0,\\
  \label{eq:asso:phi2}
  \phi(-L_2/2-d_2) = \phi(L_2/2+d_2) = 0.
\end{align}
Different from the main step of the NMM method, the separation of
variables here leads to an eigenvalue problem for $\psi$ (a function
of $x$),  instead of $\phi$ which is  not continuous at
$y=y_1$. 

The eigenvalue problem for $\psi$ can be solved by a pesudospectral
method as in \cite{sonyualu11}.  
If $M$ numerical eigenpairs 
$\{\delta_j,\psi_j(x)\}_{j=1}^M$ are obtained corresponding to the
collocation points $\{x_j\}_{j=1}^M\subset[-L_1/2-d_1,L_1/2+d_1]$, we 
approximate $\tilde{u}^s$ by
\begin{align}
  \label{eq:us:top}
  \tilde{u}^s \approx \left\{\begin{array}{lc}
                         \sum_{j=1}^M c_j^{t}e^{i\sqrt{k_0^2\varepsilon_+ + \delta_j}(y-y_1)}\psi_j(x), & y>y_1,\\
                         \sum_{j=1}^M \phi_j(y)\psi_j(x), & y_0\leq y\leq y_1,\\
                         \sum_{j=1}^M
                               c_j^{b}e^{-i\sqrt{k_0^2\varepsilon_- +
                               \delta_j}(y-y_0)}\psi_j(x), & y<y_0, 
                       \end{array}\right.
\end{align}
where the square roots have nonnegative imaginary parts. Notice that
only outgoing waves are retained in the top and bottom layers.

The functions $\phi_j(y)$ satisfy
\begin{equation}
\label{eq:phixj}
 \frac{d^2\phi_j}{dy^2} + [ k_0^2\tilde{\varepsilon}(y) + \delta_j]  \phi_j = 0.
 \end{equation}
Since $\tilde{u}^s$ satisfies the
transmission condition at $y=y_0$, we have
\[
  \phi_j(y_0) = c_j^b,\quad \phi_j'(y_0+) =
  -c_j^bi\sqrt{k_0^2\varepsilon_-+\delta_j}.
\]
Therefore, we enforce the following Robin boundary condition
\begin{equation}
  \label{eq:RBC:bot}
  \phi_j'(y_0+) = -i\sqrt{k_0^2\varepsilon_-+\delta_j}\phi_j(y_0).
\end{equation}
At $y=y_1$, we can find the coefficients $\{c_j^{\rm ps}, 
d_j^{\rm ps}\}_{j=1}^N$ such that 
\begin{align}
  \sum_{j=1}^{N}c_j^{\rm ps}\psi_j(x) &= \frac{i}{4}H_0^{(1)}(k_0n_+\rho(\tilde{x}(x),y_1)),\\
  \sum_{j=1}^{N}d_j^{\rm ps}\psi_j(x) &=
  \frac{d}{dy}\frac{i}{4}H_0^{(1)}(k_0n_+\rho(\tilde{x}(x),y))|_{y=y_1}
\end{align}
are exactly satisfied at the collocation points  $\{ x_j \}_{j=1}^N$. 
Since $u^{\rm tot}$ satisfies the transmission condition on $y=y_1$,
we have
\[
  c_j^t + c_j^{\rm ps} = \phi_j(y_1),\quad
  ic_j^t\sqrt{k_0^2\varepsilon_++\delta_j} + d_j^{\rm ps} = \phi_j'(y_1-). 
\]
Eliminating $c_j^t$, the above yields the following Robin boundary condition,
\begin{align}
  \label{eq:RBC:top}
  \phi_j'(y_1-) - i\sqrt{k_0^2\varepsilon_++\delta_j}\phi_j(y_1) = d_j^{\rm ps} - i\sqrt{k_0^2\varepsilon_++\delta_j}c_j^{\rm ps}.
\end{align}

As shown in Proposition A.2, the boundary value problem
(\ref{eq:phixj}), (\ref{eq:RBC:bot}) and  (\ref{eq:RBC:top}) has a
unique solution. Using a pesudospectral method, we solve this boundary
value problem and obtain $\phi_j(y)$ 
at the collocation points $\{y^j\}_{j=1}^N\cap[y_0,y_1]$. Finally,
since $c_j^b=\phi_j(y_0)$ and $c_j^t = \phi_j(y_1)-c_j^{\rm ps}$, 
we  have $\tilde{u}^s$ and $u_2^{\rm tot}$ in $B_2$. The other two
cases (b) and (c) are similar; we omit the details here.


\section{Numerical examples}
In this section, we carry out several numerical experiments to exhibit the
performance of our NMM method. In all examples, the physical domain is chosen to
be $(-2.5,2.5)\times(-2.5,2.5)$, and the free-space wavenumber $k_0=2\pi/\lambda$
with wavelength $\lambda=1.13$.

\begin{figure}[!ht]
  \centering
(a)\includegraphics[width=0.40\textwidth]{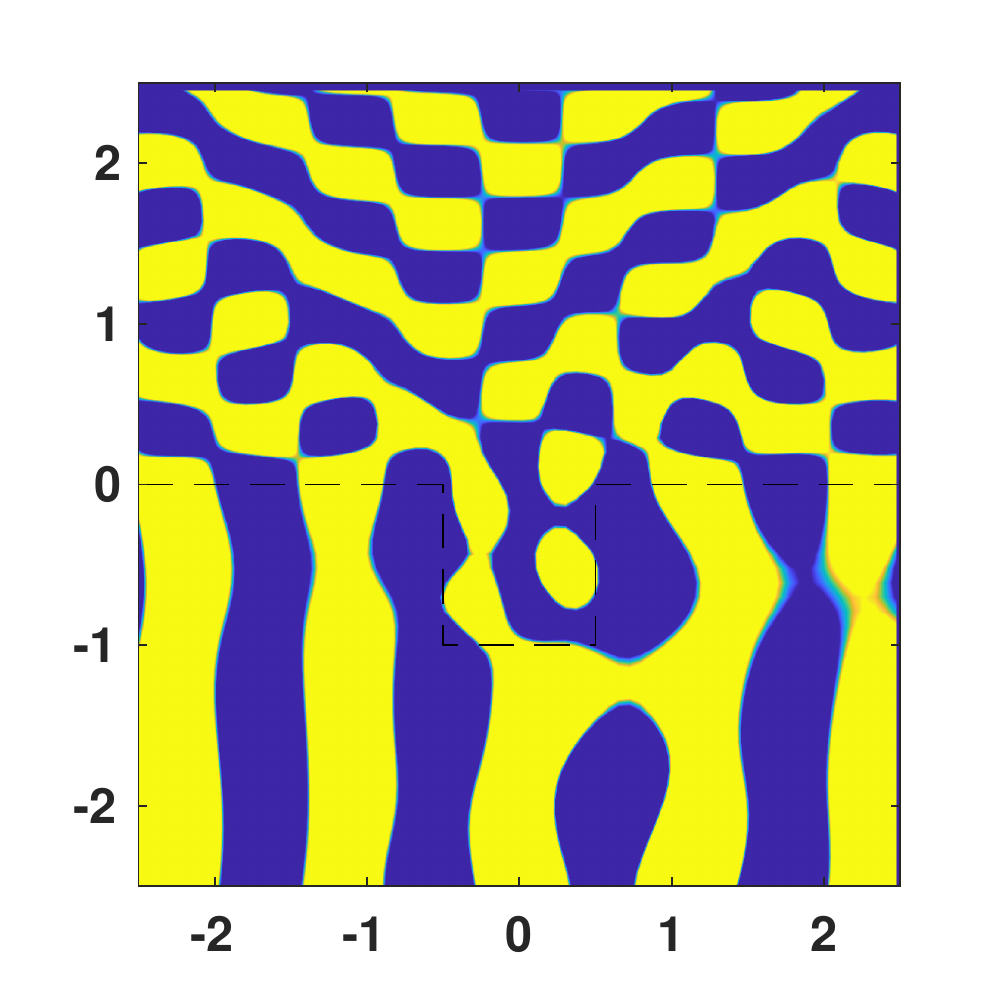}
  \includegraphics[width=0.40\textwidth]{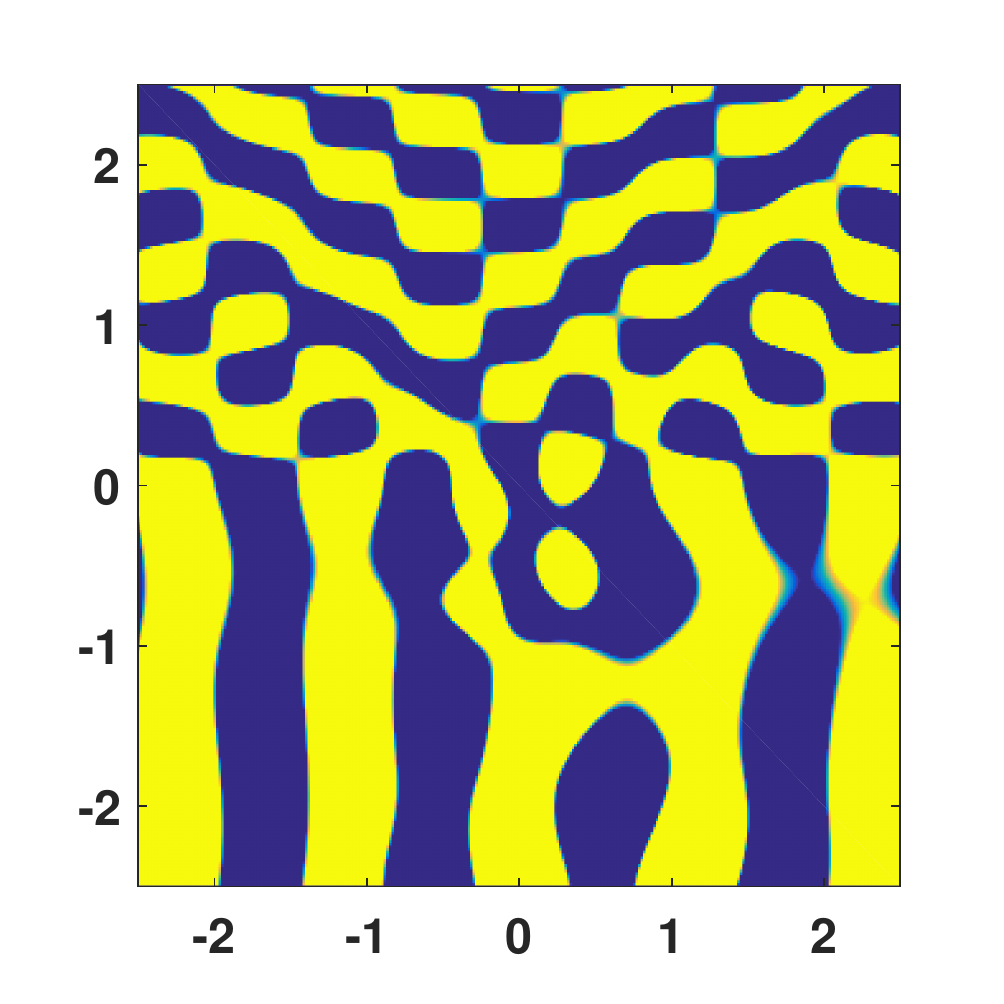}
  (b)\includegraphics[width=0.40\textwidth]{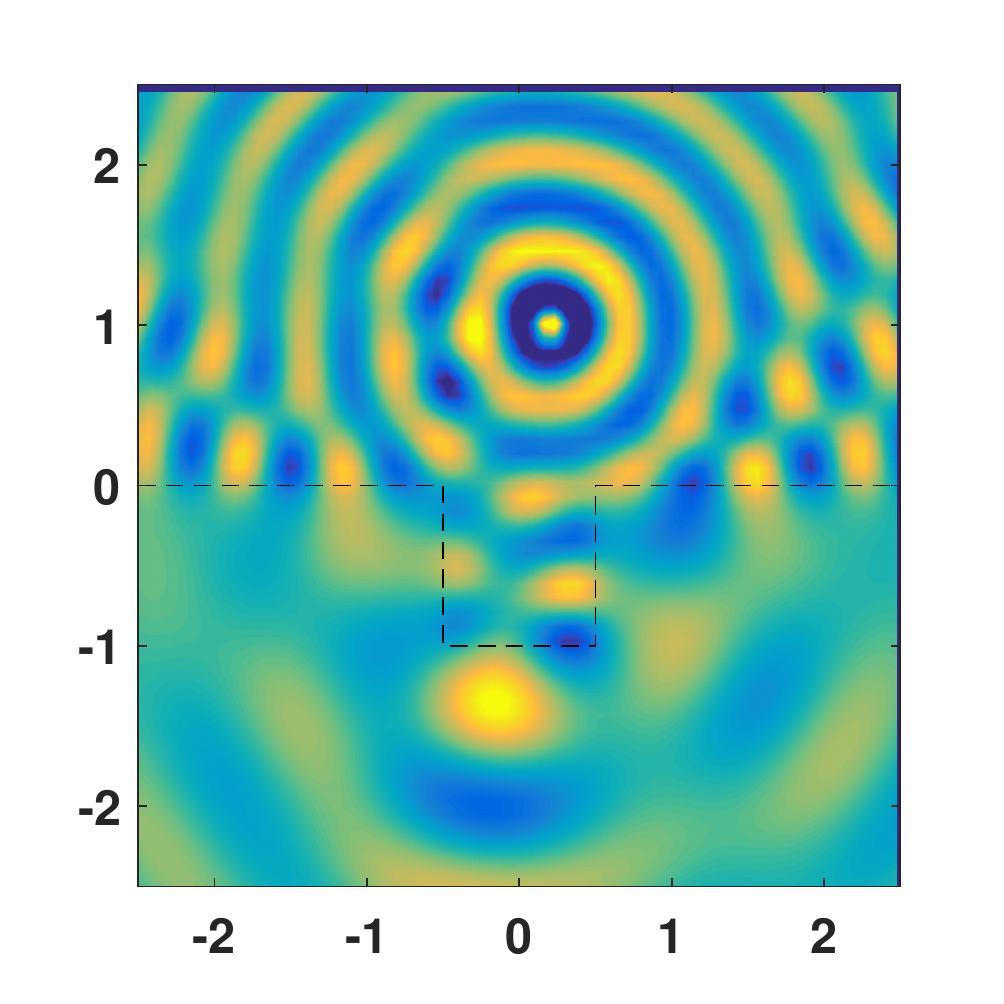}
  \includegraphics[width=0.40\textwidth]{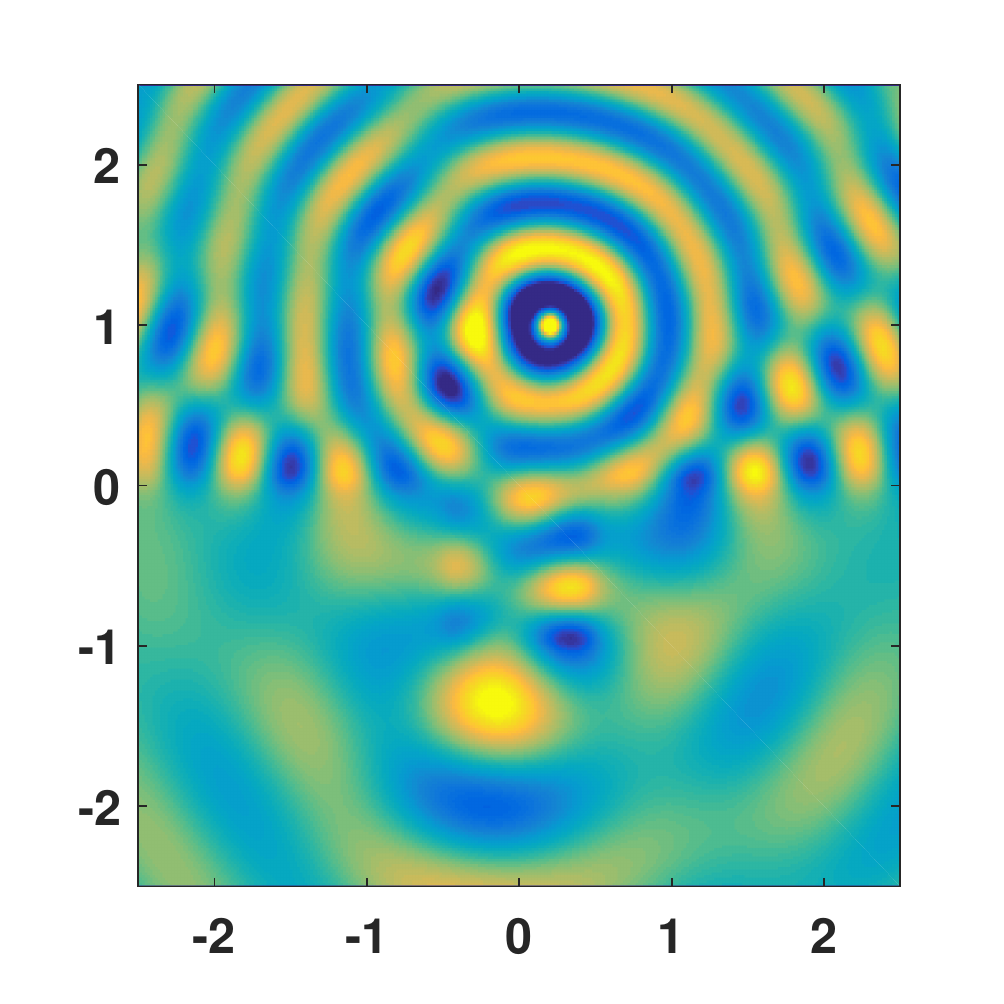}
	\caption{Example 1: Numerical solutions of a scattering
          problem in $E$ polarization: (a) plane wave at critical
          incident angle $\theta=\pi/6$; (b) cylindrical wave
          for a source at $(0.2, 1)$. For both (a) and (b), the reference solution by
          the PML-BIE method is shown on the left and the numerical
          solution by the NMM method is shown on the right.}
  \label{fig:ex1:0}
\end{figure}

\noindent{\bf Example 1.} In the first example, the background two-layer medium
is separated by interface $y=0$, with $\varepsilon(x,y)=4$ in
the top, and $\varepsilon(x,y)=1$ in the bottom. The inhomgeneity filled in
domain $D=(-0.5, 0.5)\times (-1, 1)$ is the same as the medium in the top, so
that it functions as a local perturbation to the interface $y=0$; see
the dashed lines in Fig.~\ref{fig:ex1:0}. The PML-BIE
method recently developed in \cite{luluqia18} is applicable to this problem and
is used to validate our NMM method. Using $1000$ points to discretize the
interface in the PML-BIE method, we obtain a reference solution $u^{\rm
  tot}_{\rm ref}$. To quantify the accuracy of the NMM method, we define the
following relative error,
\begin{align}
  \label{eq:rel:err}
  { e_{\rm rel}} = \frac{\max_{(x,y)\in S}|u^{\rm tot}_{\rm ref}(x,y) - u^{\rm
      tot}_{\rm NMM}(x,y)|}{\max_{(x,y)\in S}|u^{\rm tot}_{\rm ref}(x,y)|}.
\end{align}
Notice that $e_{\rm rel}$ compares the numerical solution by the NMM
method with the 
reference solution on the set $S=\{(x, y)|x=\pm 0.5, y = \pm 2.5, -1,
0\}$. The choice of $S$ is 
typical,  since it contains all corners on the interfaces and some
points at the interior boundary of the PML. 

\begin{figure}[!ht]
  \centering
  (a)\includegraphics[width=0.4\textwidth]{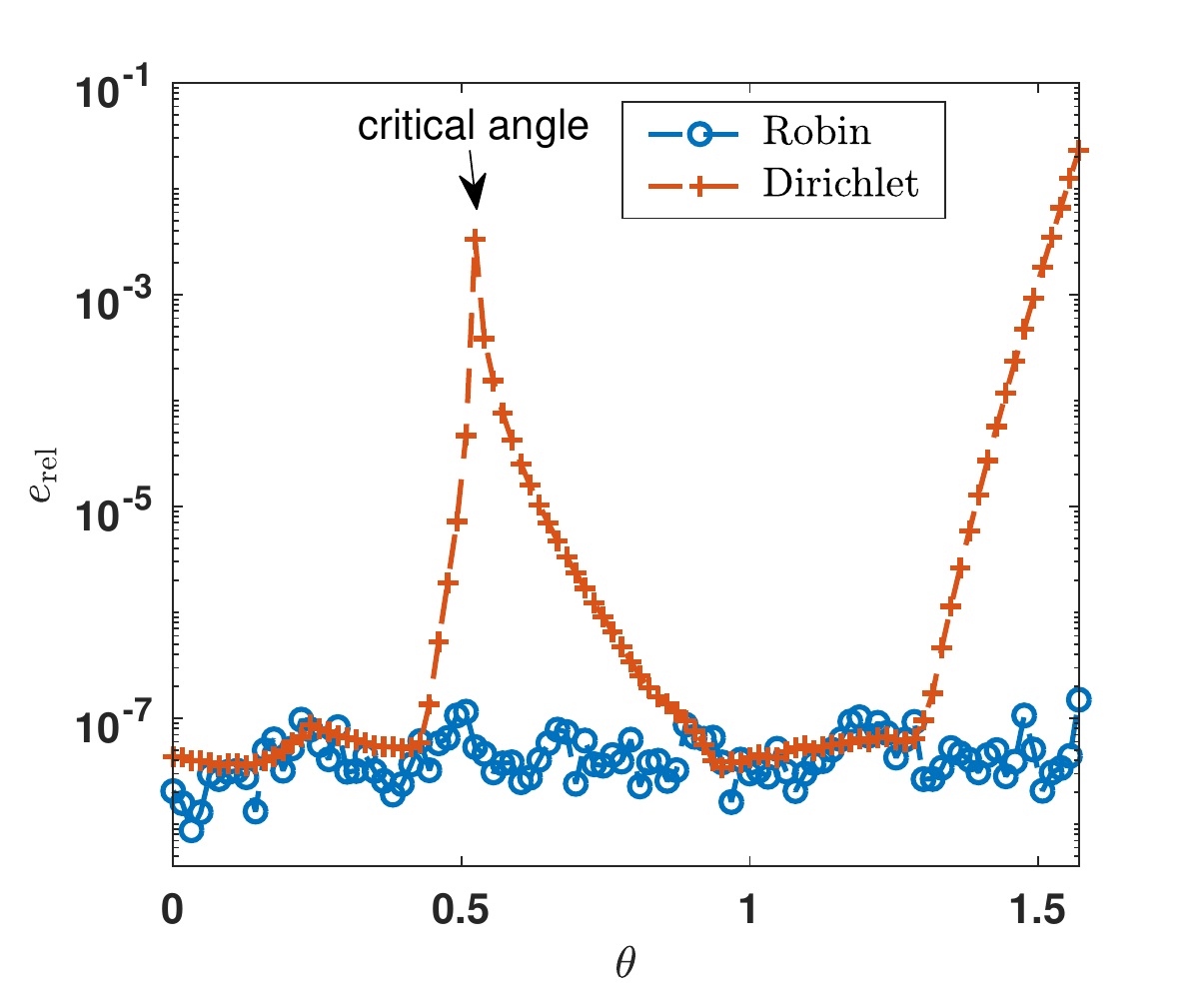}
  (b)\includegraphics[width=0.4\textwidth]{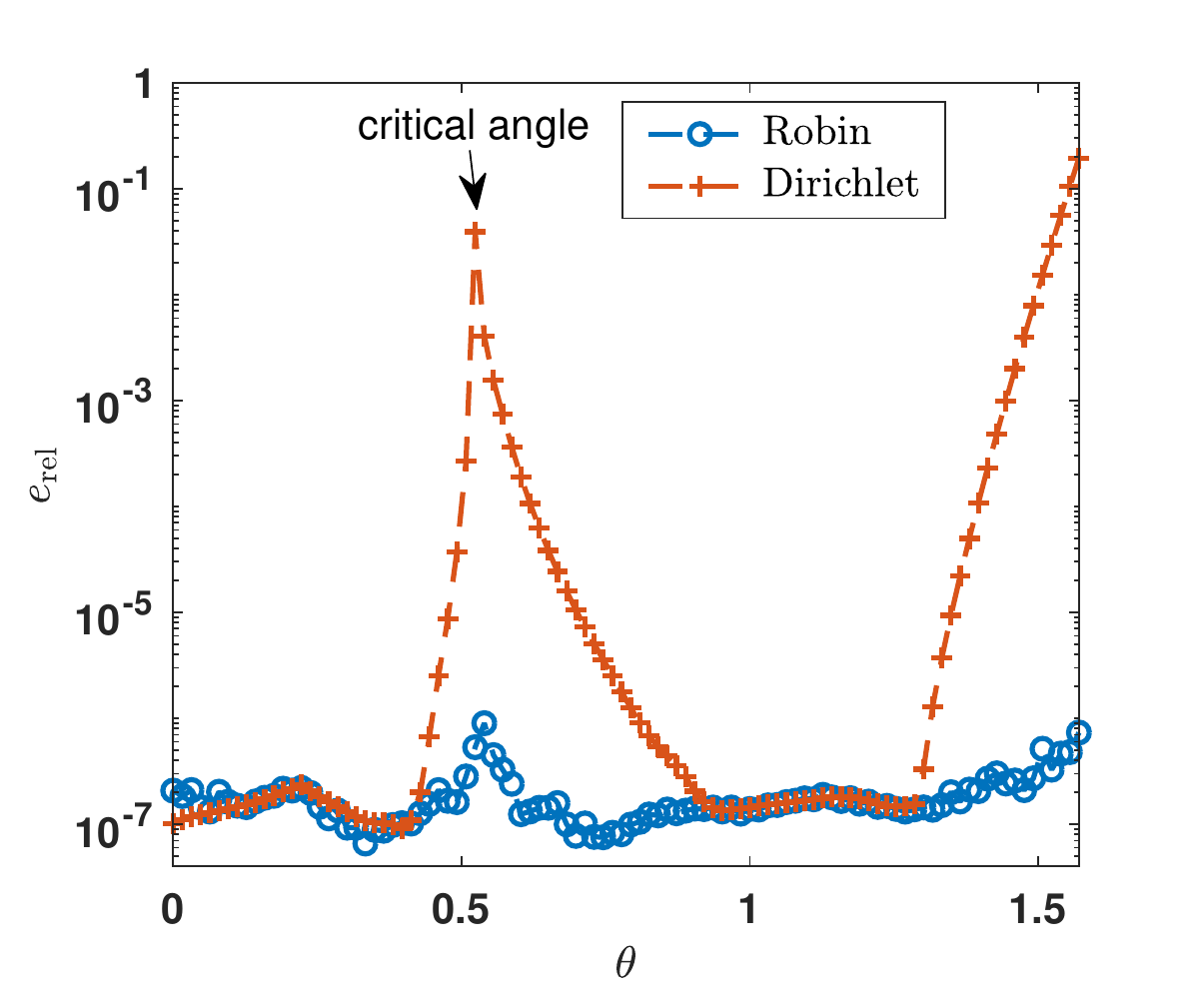}\\
  (c)\includegraphics[width=0.4\textwidth]{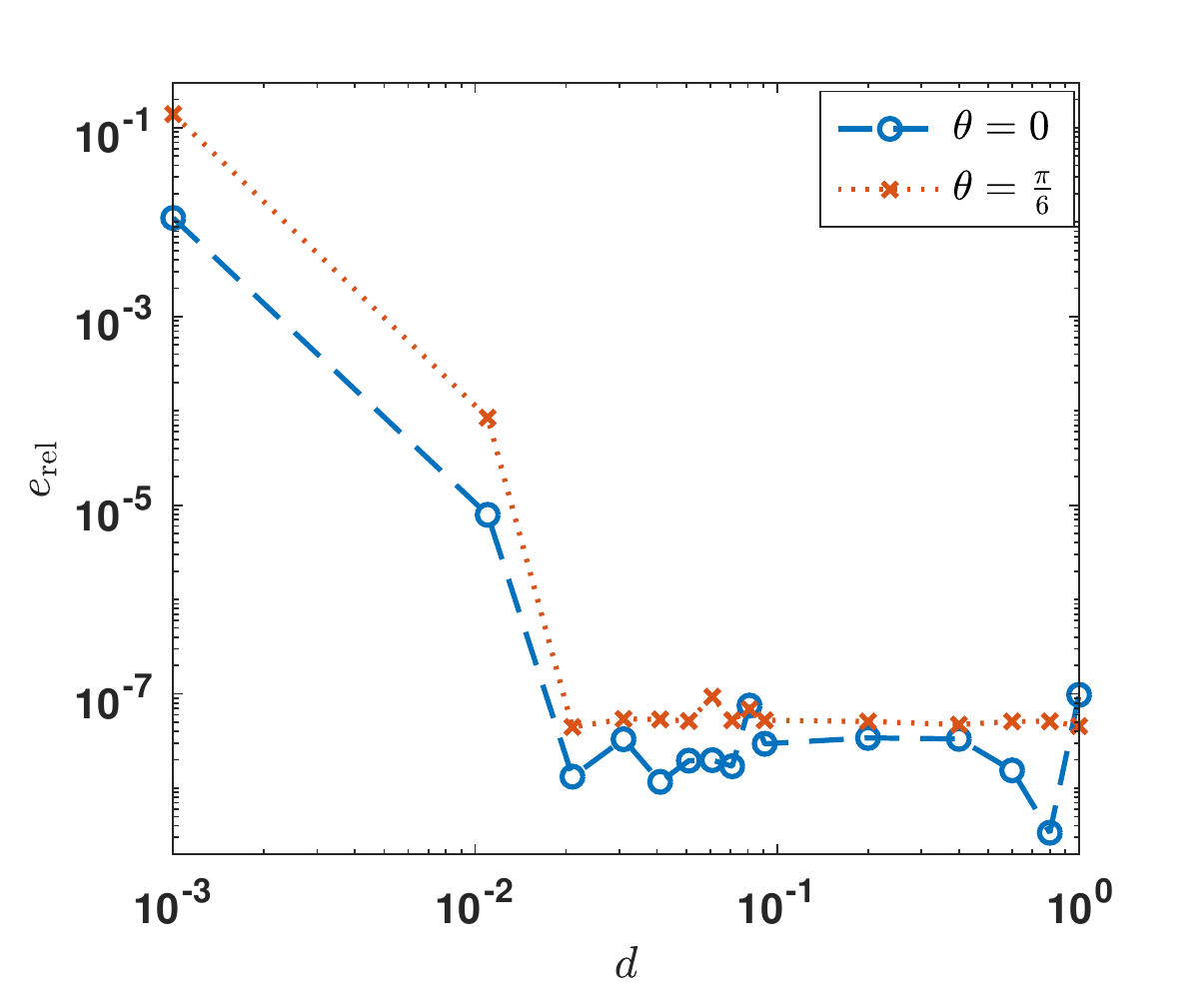}
  (d)\includegraphics[width=0.4\textwidth]{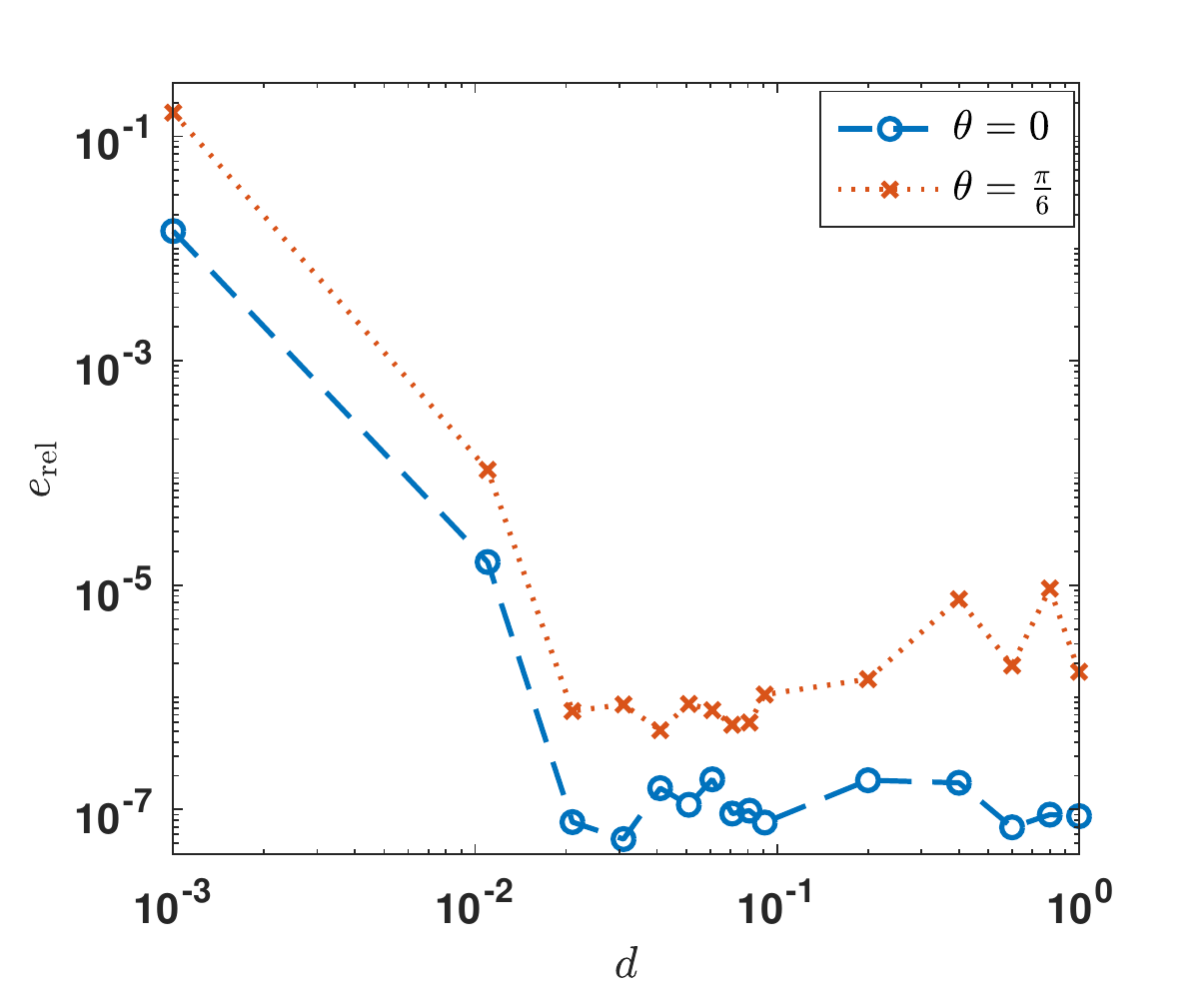}\\
	\caption{Example 1: (a) and (b): Relative error versus
          incident angle $\theta$ 
          using Robin and Dirichlet boundary conditions on
    $\Gamma_2^2$, for $N=950$, $\sigma=70$, $d=0.05$, and $m=0$.
(a) $E$ polarization;  (b) $H$     polarization. 
(c) and (d): Relative error versus PML thickness $d$ at 
          $\theta=0$ and $\theta = \pi/6$, for $N=950$, $m=0$ and 
          $\sigma=70$: (c) $E$ polarization; (d) $H$ polarization.}
  \label{fig:ex1:2}
\end{figure}

First, we validate the Robin-type boundary condition. We choose $\sigma=70$,
$d=0.05$, and $m=0$ to set up  the PML.  For both $E$ and $H$ polarizations, we
choose $N=950$ eigenmodes in each segment, and compute $e_{\rm rel}$ for 
incident angle $\theta$ varying in $[0,\pi/2)$. The results are shown in
Figs.~\ref{fig:ex1:2}(a) and \ref{fig:ex1:2}(b).
It is clear that in the vicinity of the critical angle
$\theta=\pi/6$, where total internal reflection first occurs, or
$\theta=\pi/2$, which gives horizontally propagating incident plane
waves, the Robin boundary condition produces a much smaller 
$e_{\rm rel}$ and significantly outperforms the Dirichlet boundary condition.

At both the critical incident angle $\theta=\pi/6$ and  
the normal incidence with $\theta=0$, we study the relation between 
$e_{\rm rel}$ and the PML thickness $d$ for a  
fixed $\sigma=70$. 
The numerical results are shown  
in Figs.~\ref{fig:ex1:2}(c) and \ref{fig:ex1:2}(d),
where both axes are scaled logarithmically. 
 When $d$ is small, we expect that the error is dominated by the
truncation of the PML. The results in Figs.~\ref{fig:ex1:2}(c) and 
\ref{fig:ex1:2}(d) indicate
that $e_{\rm rel}$ initially decays exponentially as $d$ is
increased. This behavior is expected from Theorem 3.1. 

Finally, we compare the numerical solutions 
by the NMM and PML-BIE methods for two types of incident waves: a
plane wave with the critical incident angle $\theta=\pi/6$,
and a cylindrical wave excited by a source at $(0.2,1)$. The results
are shown in Fig.~\ref{fig:ex1:0}. 
For each case, the PML-BIE solution is shown on the left and the NMM
solution is shown on the right. Clearly, the solutions obtained by the
two numerical methods are nearly 
indistinguishable from each other.

\begin{figure}[!ht]
  \centering
  (a)\includegraphics[width=0.38\textwidth]{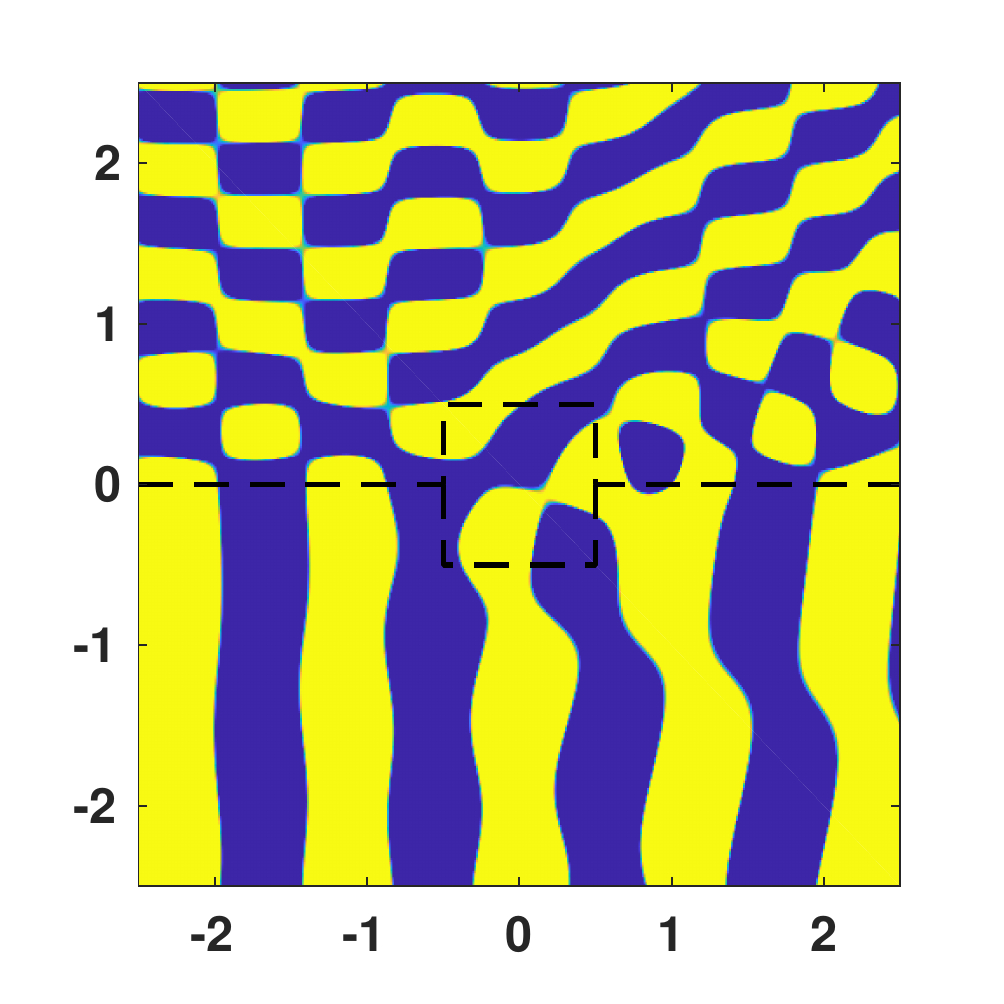}
  (b)\includegraphics[width=0.38\textwidth]{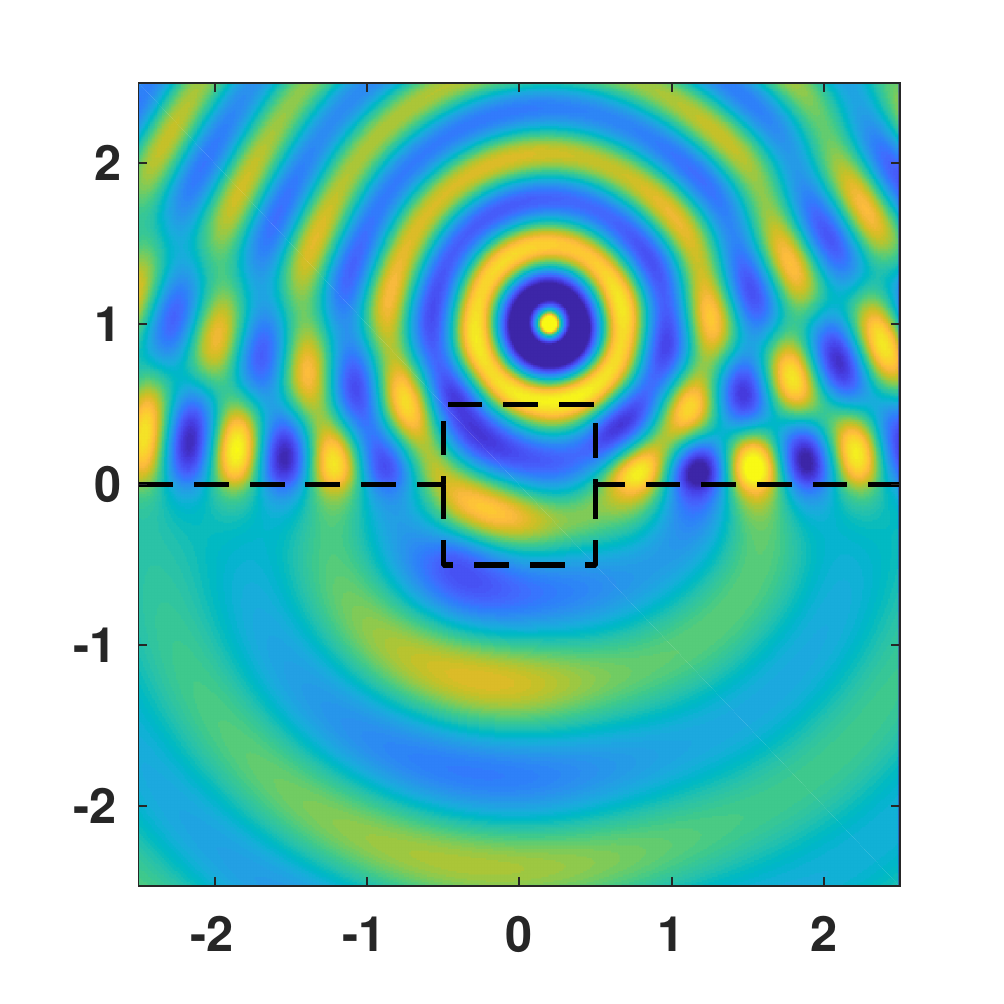}
  (c)\includegraphics[width=0.42\textwidth]{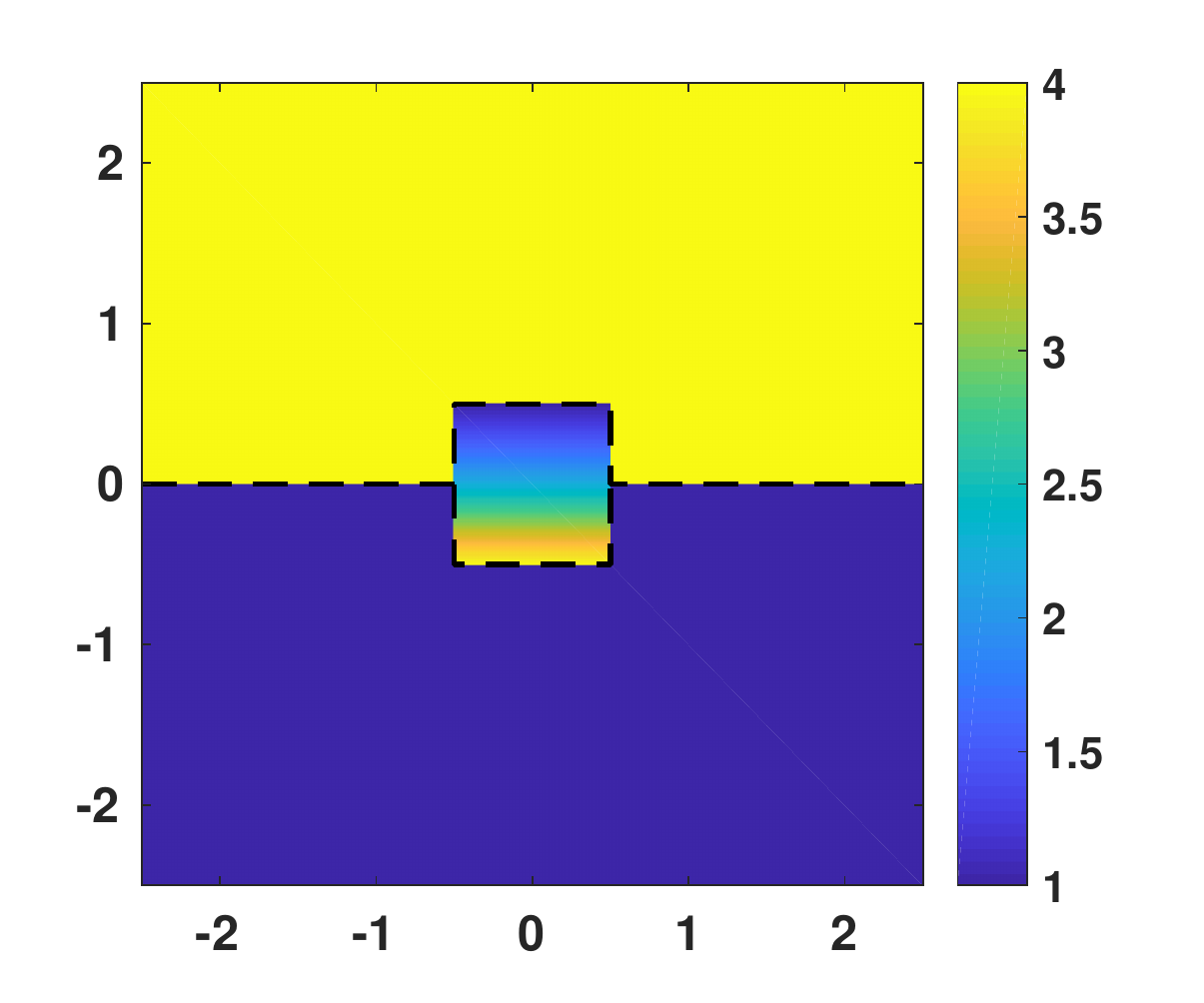}
  (d)\includegraphics[width=0.42\textwidth]{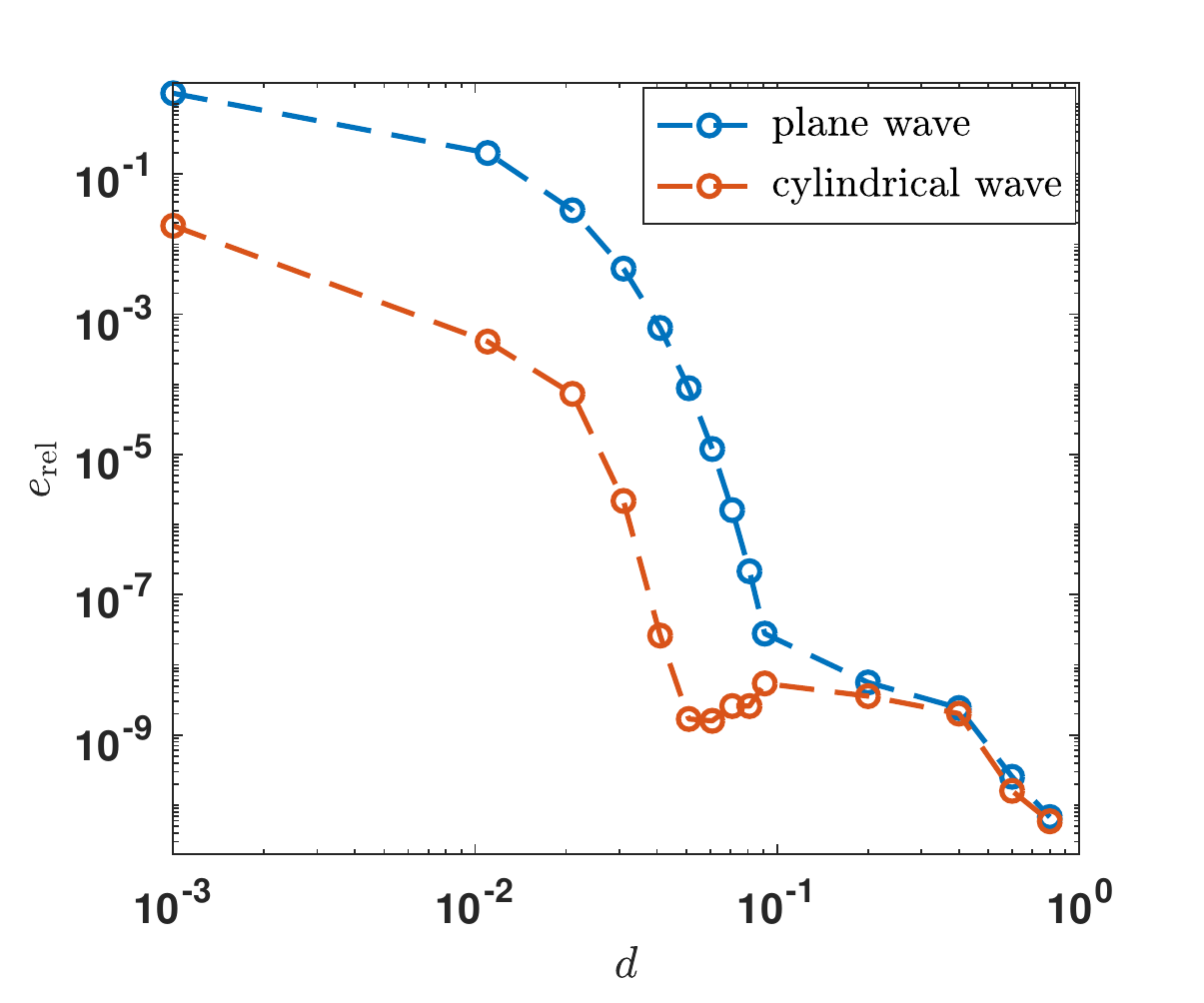}
	\caption{Example 2: Numerical solutions of a scattering
          problem in the $E$ polarization:
          (a) plane wave at critical incident angle $\theta=\pi/6$; b)
          cylindrical wave for a source at $(0.2, 1)$; (c) profile of
          the dielectric function $\varepsilon(x,y)$; (d) relative
          error $e_{rel}$ versus PML thickness $d$. 
In (a) and (b), we take $N=534$, $m=0$, $\sigma=70$ and $d=1$.}
  \label{fig:ex2:0}
\end{figure}

\noindent{\bf Example 2.} The dielectric function $\varepsilon(x,y)$ is profiled
by Fig.~\ref{fig:ex2:0}(c), where $\varepsilon(x,y)$ is $4$ in the top layer,
$1$ in the bottom layer, and $(1.5+y)^2$ in the inhomogeneity
$D=(-0.5,0.5)\times(-0.5,0.5)$.

Since $\varepsilon(x,y)$ is variable in $D$, the PML-BIE method is no longer
applicable. 
We use the NMM method to find the total field $u^{\rm tot}$ for the $E$
polarization and for two different incident waves: a plane wave at the critical
incident angle $\theta=\pi/6$, and a cylindrical wave excited by a source at
$(0.2,1)$.

For these two incident waves, using $N=534$ eigenmodes in each
segment, and using $m=0$, $\sigma=70$ and $d=1$ to set up the PML, we obtain two
numerical solutions, relatively, as shown in Figs.~\ref{fig:ex2:0}(a) and
\ref{fig:ex2:0}(b).
Using the above two numerical solutions as reference solutions, relatively, we
compute $e_{\rm rel}$ defined in Eq.~(\ref{eq:rel:err}), but for
$S=\{(x,y)|x=\pm 0.5, y=\pm 0.5, \pm 2.5\}$ for numerical solutions with values
of $d$ less than $1$, for the two incident waves, relatively; see
Fig.~\ref{fig:ex2:0}(d). As before, when $d$ is small, the relative error decays
exponentially, since it is dominated by the truncation of the PML.

\begin{figure}[!ht]
  \centering
  (a)\includegraphics[width=0.38\textwidth]{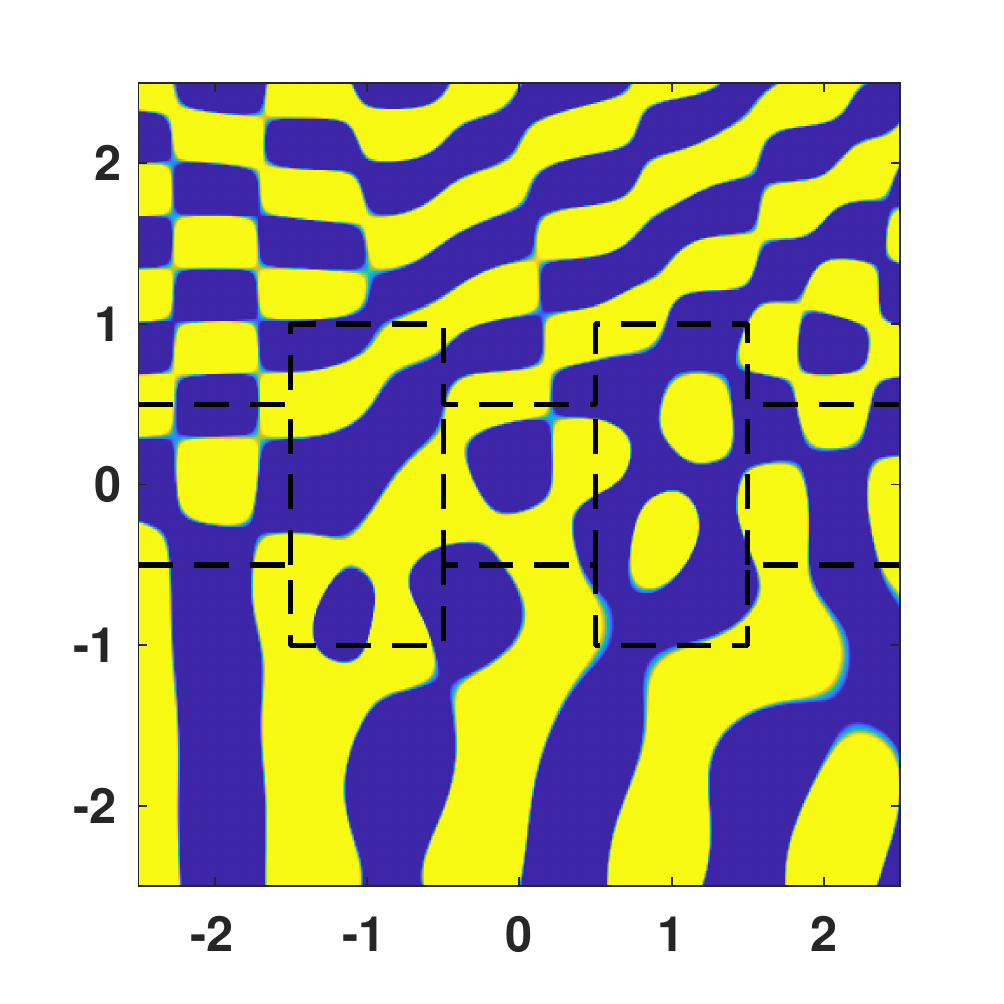}
  (b)\includegraphics[width=0.38\textwidth]{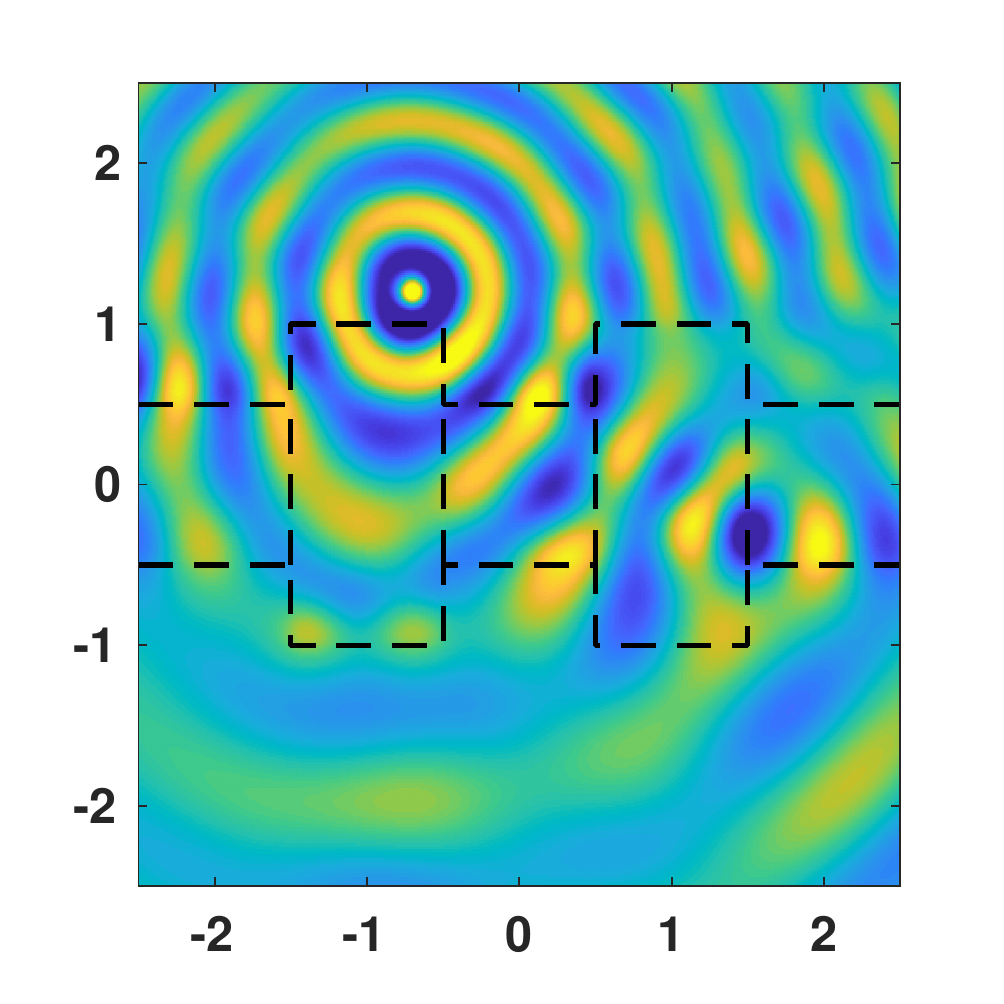}
  (c)\includegraphics[width=0.42\textwidth]{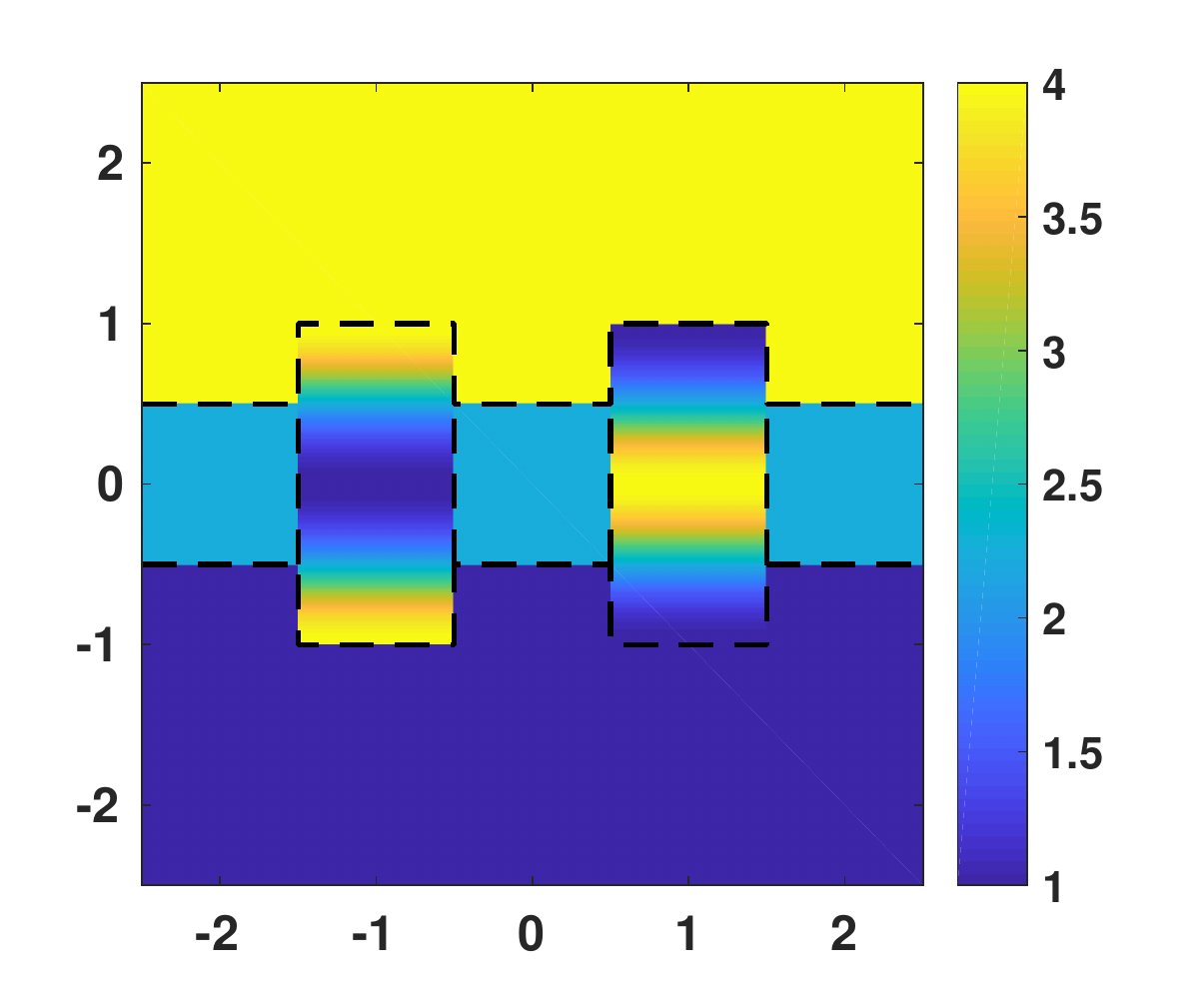}
	\caption{Example 3: Numerical solutions of a scattering
          problem in the $E$ polarization: (a) plane wave at critical
          incident angle $\theta=\pi/6$; (b) cylindrical wave
          for a source at $(-0.7, 1.2)$; (c) profile of the dielectric
          function $\varepsilon(x,y)$. Other parameters: $N=633$, $m=0$,
          $\sigma=70$, and $d=1$.}
  \label{fig:ex3:0}
\end{figure}

\noindent{\bf Example 3.} The dielectric function $\varepsilon(x,y)$ is profiled
by Fig.~\ref{fig:ex3:0} (c), where two $y$-dependent inhomogeneities are embedded in
the background medium with three layers. Here, $\varepsilon(x,y)$ is $4$, $2.25$
and $1$ in the top, inner, and bottom layers, respectively, and in the two
inhomogeneities,
\[
  \varepsilon(x,y)=\left\{\begin{array}{ll}
(1+\sin^2(\pi y/2))^2, & (x,y)\in D_1=(-1.5,-0.5)\times(-1,1);\\
(1+\cos^2(\pi y/2))^2, & (x,y)\in D_2=(0.5,1.5)\times(-1,1).\\
                            \end{array}\right.
\]

Using $N=633$ eigenmodes in each segment, and using $m=0$, $\sigma=70$, and
$d=1$ to set up the PML, we calculate the total field $u^{\rm tot}$ for the $E$
polarization and for two different incident waves: a plane wave with the
critical incident angle $\theta= \pi/6$, and a cylindrical wave excited by a
source at $(-0.7, 1.2)$. The results are shown in Figs.~\ref{fig:ex3:0}(a) and
\ref{fig:ex3:0}(b).

\section{Conclusion}

The NMM methods are widely used in engineering applications for simulating
propagation and scattering of linear electromagnetic, acoustic and elastic
waves. These methods are restricted to special structures, but are more
efficient than the standard numerical methods when they are applicable, since no
discretization is needed for one spatial variable. In this paper, a new NMM
method is developed to overcome a limitation of existing NMM methods due to the
existence of a non-propagating and non-decaying wave field component. A
Robin-type boundary condition is used to ensure that the wave field component
with zero or near-zero transverse wavenumber is not reflected by a truncated
PML. A theoretical foundation of the new NMM is established by a theorem which
reveals the exponential convergence of the PML solution with the hybrid
Dirichlet-Robin boundary conditions. In addition, for scattering problems with
cylindrical incident waves, we developed a fast method to compute reference
solutions needed in the NMM methods. Numerical examples are presented to
validate the NMM method and illustrate its performance.

We have implemented the new NMM method for two-dimensional structures with one
or more inhomogeneities, for electromagnetic waves in $E$ and $H$ polarizations,
and for both plane and cylindrical incident waves. The NMM methods are also
applicable to three-dimensional (3D) rotationally symmetric structures that are
piecewise uniform in the radial variable \cite{hulinluoh18, lushilu14}. There is
also a related method for more general 3D structures without the rotational
symmetry \cite{shilu15, shilulu16}. The techniques developed in this paper,
namely, the Robin-type boundary condition for terminating the PML and the fast
method for computing reference solutions for cylindrical incident waves, should
also be useful in these NMM and related methods for 3D structures.


\begin{appendices}


\section{Proof of Propositions}

\begin{prop}
  Under the same assumptions as Theorem 2.2, we consider the eigenvalue problem
  (\ref{eq:phi:te}-\ref{eq:phi:te3}) with $\varepsilon_i$ replaced by
\[
  \varepsilon_{\rm gen}(y) = \left\{
    \begin{array}{lc}
      \varepsilon_+,& y>y_1,\\
      \varepsilon_-,& y<y_0,\\
      \varepsilon_{\rm PHY}(y), & y_0\leq y\leq y_1, 
    \end{array}
  \right.
\]
where $\varepsilon_{\rm PHY}(y)$ can be any piecewise smooth and positive
function. Then, one and only one of the following two cases occurs:
  \begin{itemize}
  \item[(a)] There exist $\sigma^0>0$ and $d^0>0$ such that if
    $\sigma>\sigma^0$ or if $d>d^0$, then ${\rm Im}(\delta)\geq 0$ for any
    eigenpair $\{\phi,\delta\}$ that solves (\ref{eq:phi:te}-\ref{eq:phi:te3});
  \item[(b)] For a fixed $d>0$ ($\sigma>0$), there exist a sequence of
    $\{\sigma^n\}_{n=1}^{\infty}$($\{d^n\}_{n=1}^{\infty}$, respectively) that
    approaches infinity as $n\rightarrow\infty$ such that there exists a sequence
    of associated eigenpairs $\{\phi^n,\delta^n\}$ satisfying ${\rm
      Im}(\delta^n)<0$.
  \end{itemize}
  If case (b) holds, then ${\rm Im}(\delta^n)\rightarrow 0$ as
  $n\rightarrow\infty$ and for sufficiently large $n$,
  \[
    {\rm Re}(\delta^n)\in[\max(k_0^2\varepsilon_+,k_0^2\varepsilon_-),
    \max(k_0^2\varepsilon_{\rm PHY})).
  \]
  \begin{proof}
    It is clear that if case (a) does not hold, then case (b) must hold. We now
    prove ${\rm Im (\delta^n)}\rightarrow 0$ as $n\rightarrow\infty$. Integrating
    (\ref{eq:phi:te}) with $\bar{\phi}$ on $[-L_2/2,L_2/2]$ and using
    integration by parts yield
    \begin{align}
      \label{eq:det:delta}
    \delta\int_{-L_2/2}^{L_2/2} |\phi|^2dy =& -\int_{-L_2/2}^{L_2/2}|\phi'(y)|^2dy +
      k_0^2\int_{-L_2/2}^{L_2/2}\varepsilon_{\rm PHY}|\phi|^2dy \nonumber\\
      &+ \left.\left(
        \frac{d\phi}{dy}\bar{\phi} \right)\right|_{-L_2/2}^{L_2/2}.
      \end{align}
      In ${\rm PML}_y=(-L_2/2-d,-L_2/2)\cup(L_2/2,L_2/2+d)$, $\phi(y)$ has the
      following general solution form,
   \[
    \phi(y) = c^1_{\pm}e^{\pm ik^*_{\pm} ( \tilde{y}\mp L_2/2)} +
    c^2_{\pm}e^{\mp ik^*_{\pm} ( \tilde{y}\mp L_2/2)},\quad{\rm in}\ {\rm PML}_{y}\cap \mathbb{R}_{\pm},
  \] 
  where $k^*_{\pm}=\sqrt{k_0^2\varepsilon_{\pm}-\delta}$ with ${\rm
    Re}(k^*_{\pm})\geq 0$. The homogeneous Dirichlet boundary condition at
  $y=\pm(L_2/2+d)$ implies that
  \[
    c^2_{\pm} = -c^1_{\pm} e^{\pm i 2k_{\pm}^* d(1+i\sigma)},
  \]
  so that, by a straightforward calculation, one obtains
  \begin{align}
    \label{eq:est:imag}
    \left.\left(
    \frac{d\phi}{dy}\bar{\phi} \right)\right|_{-L_2/2}^{L_2/2} = \sum_{l=\pm}|c_l^1|^2\Big[(1 - e^{-4d({\rm Im}(k^*_{l}) + {\rm Re}(k^*_l)\sigma)})ik^*_l \nonumber\\
    - 2k^*_l e^{-2d({\rm Im}(k^*_l) + {\rm Re}(k^*_l)\sigma)}\sin(2d({\rm Re}(k^*_l)-{\rm Im}(k^*_l)\sigma))\Big].
  \end{align}
  If ${\rm Im}(\delta)<0$, then ${\rm Im}(k_{\pm}^{*})>0$. Next, we show that
  ${\rm Re}(k_\pm^{*,n})\rightarrow 0$ as $n\rightarrow\infty$, where
  $k_{\pm}^{*,n}=\sqrt{k_0^2\varepsilon_{\pm}-\delta^n}$ with nonnegative real
  part. Otherwise, suppose we have a subsequence $\{n_k\}_{k=1}^{\infty}$ such
  that ${\rm Re}(k_{\pm}^{*,n_k})\rightarrow c_{\pm}^{0}> 0$ as
  $n_k\rightarrow\infty$. Considering the imaginary part of (\ref{eq:est:imag}),
  \[
 {\rm Im}\left.\left(
    \frac{d\phi^{n_k}}{dy}\bar{\phi}^{n_k} \right)\right|_{-L_2/2}^{L_2/2}\rightarrow
|c_+^1|^2c_+^0 + |c_-^1|^2c_-^0> 0,    
  \]
  which implies that ${\rm Im}(\delta^{n_k})>0$. Considering the real part of
  (\ref{eq:est:imag}),
  \[
{\rm Re}\left.\left(\frac{d\phi}{dy}\bar{\phi}
  \right)\right|_{-L_2/2}^{L_2/2}\leq\sum_{l=\pm}|c_l^1|^2\Big[-(1 -
e^{-4d({\rm Im}(k^*_{l}) + {\rm Re}(k^*_l)\sigma)}){\rm Im}(k_l^*) \nonumber\\
+ 2{\rm Re}(k^*_l)\Big].
\]
For sufficiently large $n$, ${\rm Re}(k^{*,n}_l)$ can be arbitrarily small such
that
\[
{\rm Re}\left.\left(\frac{d\phi^{n}}{dy}\bar{\phi}^n
  \right)\right|_{-L_2/2}^{L_2/2}\leq 0.
\]
Therefore, considering the real part of (\ref{eq:det:delta}),
\[
  {\rm Re}(\delta^n)< k_0^2\max(\varepsilon_{\rm PHY}),
\]
since $\frac{d}{dy}\phi^{n}(y)\neq 0$ in $[-L_2/2,L_2/2]$. On the other hand,
  \[
    {\rm Re}(\delta) = k_0^2\varepsilon_{\pm} - {\rm Re}(k^*_{\pm})^2 + {\rm
      Im}(k^*_{\pm})^2\geq k_0^2\varepsilon_{\pm} - {\rm Re}(k^*_{\pm})^2.
  \]
  so that ${\rm Re}(\delta^n)\geq k_0^2\varepsilon_{\pm}$ and that
  \[
    {\rm Im}(k_{\pm}^{*,n})<k_0^2\max(\varepsilon_{PHY})-k_0^2\varepsilon_{\pm} +
    {\rm Re}(k_{\pm}^{*,n})^2,\]
  Consequently, as $n\rightarrow \infty$, we have
  ${\rm Im}(\delta^n)=-2{\rm Im}(k_{\pm}^{*,n}){\rm
    Re}(k_{\pm}^{*,n})\rightarrow 0$.
  \end{proof}
\end{prop}
\begin{prop}
  If $\tilde{ \varepsilon }(y)$ is smooth on $[y_0,y_1]$,  the following field
 \begin{align}
   \label{eq:u2:closed:form}
   u^{\rm tot}_2 = e^{-i\beta_+ y_1}\left\{
   \begin{array}{ll}
     e^{i\alpha x}( e^{-i\beta_+ (y-y_1))} + R_2e^{i\beta_+ (y-y_1)})& {\rm if}\ y\geq y_1,\\
e^{i\alpha x} f(y),&{\rm if}\ y_0<y<y_1,\\
     T_2e^{i(\alpha x - \beta_-(y-y_0))},&{\rm if}\ y\leq y_0.
   \end{array}
   \right.
 \end{align}
 solves the scattering problem (\ref{eq:gov:problem}) and (\ref{eq:trans:cond})
 with $\varepsilon(x,y)=\varepsilon_2(y)$ in $\mathbb{R}^2$, where $\beta_-$ was
 defined in (\ref{eq:sol:para}), the unknown function $f\in C^{2}[y_0,y_1]$ is
 uniquely determined by the following boundary value problem
 \begin{align}
   \label{eq:bvp:1}
   &f'' + (k_0^2\varepsilon_2  - \alpha^2) f = 0,\\
   \label{eq:bvp:2}
   &f'(y_0+) = -i\beta_-f(y_0), \\
   \label{eq:bvp:3}
   &f'(y_1-) = i\beta_+(f(y_1)-2),
 \end{align}
 and 
 \[
   R_2 = f(y_1) - 1,\quad T_2 = f(y_0).
 \]
 \begin{proof}
   The verification that $u_2^{\rm tot}$ defined in (\ref{eq:u2:closed:form}) is
   indeed a solution is straightforward. One only needs to prove that the
   boundary value problem (\ref{eq:bvp:1}-\ref{eq:bvp:3}) has a unique solution.
   By the standard ODE theory, one needs to show that equation (\ref{eq:bvp:1})
   with the following homogeneous Robin boundary conditions
   \[
     f'(y_0+) = -i\beta_-f(y_0),\quad f'(y_1-) =
     i\beta_+f(y_1),
   \]
   has only the trivial solution $f=0$. To show this, integrating
   (\ref{eq:bvp:1}) with $\bar{f}$ on $[y_0,y_1]$ yields, by integration by
   parts,
   \begin{align*}
     \int_{y_0}^{y_1} f'\bar{f}'dy -
     (k_0^2\tilde{\varepsilon}-\alpha^2)f\bar{f}dy -
     (i\beta_+f(y_1)\bar{f}(y_1) + i\beta_-f(y_0)\bar{f}(y_0)) = 0.
   \end{align*}
   Considering the imaginary part of the left-hand side, we have
   \[
     \beta_+|f(y_1)|^2 + \beta_-|f(y_0)|^2 = 0,
   \]
   so that $f(y_1)=0$ since $\beta_+>0$. Notice that $\beta_-$ could be zero
   when total internal reflection occurs. Then, $f'(y_1-)=i\beta_+f(y_1)=0$
   indicates that $f=0$ on $[y_0,y_1]$ which completes the proof. 
 \end{proof}
\begin{myremark}
  $u_2^{\rm tot}$ in (\ref{eq:u2:closed:form}) plus any guided mode or any surface
  mode, if there exists, still solves the scattering problem.
\end{myremark}
\end{prop}

\end{appendices}

\bibliographystyle{plain}
\bibliography{wt}
\end{document}